\documentclass[12pt]{article}
\usepackage{amsmath,amsfonts,epsfig}
\title{The Even Isomorphism Theorem for Coxeter Groups}
\author{M. Mihalik}
\newtheorem{theorem}{Theorem}
\newtheorem{proposition}[theorem]{Proposition}
\newtheorem{lemma}[theorem]{Lemma}

\date{July 16, 2003}
\begin{document}
\maketitle

\section{Introduction.}

A {\it Coxeter system} is a pair $(W,S)$ such that $W$ is a group
with {\it Coxeter presentation} $\langle S : (st)^{m_{st}}\forall
s,t\in S \rangle$ where for all $s,t\in S$, $m_{st}\in
\{1,2,\ldots ,\infty \}$, $m_{st}=m_{ts}$ and $m_{st}=1$ if and
only if $s=t$. (The relation $(st)^{\infty }$ means that $st$ has
infinite order in $W$.) Note that all $s\in S$ are order 2 and
that if $m_{st}=2$ then $s$ and $t$ commute.

There are two diagrams associated to a Coxeter group that appear
regularly in the literature. The diagrams $V_D(W,S)$ and
$V_F(W,S)$ for the Coxeter system $(W,S)$ are labeled graphs with
vertex set $S$. In $V_F$ there is an edge labeled $m_{st}$
between distinct vertices $s$ and $t$ if and only if $m_{st}\not
=\infty $. In $V_D$ there is a edge labeled $m_{st}$ between
distinct vertices $s$ and $t$ if and only if $m_{st}\ne 2$. The
vertices of components of $V_D$ generate factors of a direct
product decomposition of $W$ and the vertices of components of
$V_F$ generate factors of a free product decomposition of $W$.
While it is traditional to call $V_D$ the Coxeter graph or
Coxeter diagram for $(W,S)$, in this paper we only consider $V_F$
diagrams and we call such diagrams, Coxeter diagrams or simply
diagrams.

A Coxeter presentation is {\it even} if all $m_{st}$ for $s\not=t$
are even or $\infty$. In this case we call the corresponding
Coxeter group and diagram even.

A Coxeter group is {\it rigid} if any two Coxeter presentations
for this group are isomorphic presentations. In \cite{B}, Bahls
shows that any Coxeter group can have at most one even
presentation. In $\S$ 4, we classify the even rigid Coxeter
groups. Rigidity and a variety of analogous notions are designed
to give insight into a fundamental problem in the theory of
Coxeter groups.

\medskip

{\bf The Coxeter Isomorphism Question:} Given two Coxeter
presentations, do they present isomorphic groups?

\medskip

The Coxeter presentations $\langle
x,y,z:x^2,y^2,z^2,(xy)^3,(xz)^2,(yz)^2\rangle$ and $\langle a,b:
a^2, b^2,(ab)^6\rangle$ present isomorphic groups, but only the
latter is even. In particular, the Coxeter group presented here
is not rigid.

In this paper we produce an algorithm to decide if an arbitrary
Coxeter presentation presents a finitely generated even Coxeter
group. Furthermore, we can decide if two (finite) Coxeter
presentations present the same even Coxeter group. This solves
the even Coxeter isomorphism question.

Our Proposition \ref{P7} is used by Patrick Bahl's in his thesis
\cite{B} to show that there is an unique even Coxeter presentation
for a finitely generated even Coxeter group. We in turn use
Bahl's result in the final stage of our algorithm to decide if
two finite Coxeter presentations present the same even Coxeter
group. The proof of a vital combinatorial lemma hinges on the
visual decomposition theorem of \cite{MT}.  A critical tool in our
algorithm is that of twisting in Coxeter diagrams. This method of
producing different Coxeter diagrams (and different
presentations) for the same Coxeter group was introduced by N.
Brady, J. McCammond, B. Muhlherr and W. Neumann in \cite{BMMN}. At
this time the only known way to produce different Coxeter
diagrams for a given Coxeter group is by twisting or by a certain
triangle/edge exchange.

Our main theorem is the following:

\begin{theorem} \label{T1}
Suppose $(W,S)$ is an even Coxeter system and $V'$ is a Coxeter
diagram for $W$ with odd labeled edge $[xy]$. Then there is a
diagram $V''$ for $W$ obtained from $V'$ by first performing a
twist around $[xy]$ and then replacing a triangle $[xyu]$ by an
edge with even label.
\end{theorem}

The proof of our theorem specifically defines the set to be
twisted and a vertex $u$ so that triangle $[xyu]$ may be replaced
by an even edge. The resulting diagram for $W$ has (one) fewer odd
labeled edges than the original. Thus we have a simple algorithm
to change a non-even diagram for a finitely generated even Coxeter
group $W$ into the unique even diagram for $W$.

If $(W,S)$ is an arbitrary Coxeter system with diagram $V$
containing an odd edge $[xy]$, then either the described twist
and triangle replacement can be carried out or $W$ is not an even
Coxeter group. Hence one can decide if a given finitely generated
Coxeter group is even or not. Given two finitely generated
Coxeter systems $(W_1,S_1)$ and $(W_2,S_2)$ with diagrams $V_1$
and $V_2$ respectively, one can decide if $W_1$ and $W_2$ are
isomorphic even Coxeter groups. Simply apply our algorithm
repeatedly to $V_1$ and $V_2$ until either an odd edge cannot be
replaced by an even one using our technique (in which case one of
the groups is not even), or until all odd edges are replaced in
both diagrams. In the later scenario, Bahls' even rigidity result
implies that the resulting even diagrams are diagram isomorphic
if and only if $W_1$ and $W_2$ are isomorphic.

It is also evident that given a Coxeter system for a finitely
generated even Coxeter group, one can use the methods of this
paper to produce all other Coxeter systems for that group.

\section{Preliminaries.}

In this section, we describe: twisting in Coxeter diagrams as
introduced in \cite{BMMN}, visual decompositions of Coxeter
groups \cite{MT}, and techniques to construct quotient maps of
Coxeter groups that match quotient maps of Coxeter diagrams.

\medskip

\noindent {\bf (1) Twisting.} In an arbitrary Coxeter system
$(W,S)$, twisting makes sense for any subset of $S$ that
generates a finite subgroup of $W$. We only need twist around
pairs of distinct vertices $x,y\in S$ such that $m_{xy}$ is an
odd integer.

Suppose $V$ is a Coxeter diagram for the Coxeter system $(W,S)$.
Given $x,y\in S$, define $lk(x)$ (the ``link" of $x$) to be the
set of all vertices of $V$ that are connected to $x$ by an edge.
Define $lk_2(x)$ (the ``2-link" of $x$) to be the set of all
vertices of $V$ that are connected to $x$ by an edge labeled 2.
Define $st(x)$ (the ``star" of $x$) to be $lk(x)\cup \{x\}$.
Define $lk_2(x,y)$ (the ``2-link" of $x$ and $y$) to be
$lk_2(x)\cap lk_2(y)$, i.e. the set of all vertices in $V$ that
are connected to both $x$ and $y$ by an edge labeled 2. So each
$s\in lk_2(x,y)$ commutes with both $x$ and $y$. Denote
$b^{-1}ab$ by $a^b$. Now suppose $x$ and $y$ are distinct
elements of $S$ and $m_{xy}=2n+1$. Let $d$ be the (unique)
element of length $2n+1$ in $\langle x,y\rangle$. Note that
$x^{d}=y$ and $y^d=x$. Suppose $U\subset S-\{x,y\}$ and for each
edge $[us]$ such that $s\in S-(U\cup \{x,y\})$ and $u\in U$,
$s\in lk_2(x,y)$. Then the twisting theorem of \cite {BMMN}
implies that $(W,S')$ is a Coxeter system, where $S'=U^d\cup
(S-U)$ and a diagram for $(W,S')$ is obtained from $V$ by
changing each edge of $V$ that connects to a vertex $u\in U$ to a
vertex $v\in \{x,y\}$ to connect instead from $u$ to $v^d$, and
leaving other edges unchanged.

\medskip

\noindent {\bf (2) Visual Decompositions of Coxeter Groups.}
Suppose $V$ is a diagram for the Coxeter system $(W,S)$ and some
subset $C$ of $S$ separates vertices of $V$, then a simple
examination of presentations, shows that $W$ decomposes as
$\langle A\rangle \ast_{\langle C\rangle}B$, where $A\cup B=S$,
$A$ is $C$ union the vertices of some set of components of $V-C$
and $B$ is $C\cup (S-A)$. This type of decomposition extends in a
natural way to graphs of groups decompositions of $W$. Such
decompositions are called ``visual" decompositions of $W$ since
they are easily seen in $V$ and the main theorem of \cite{MT}
states that given any graph of groups decomposition of $W$ there
is a visual decomposition that basically refines the given
decomposition. More specifically, any vertex (edge) group of the
visual decomposition of $W$ is a subgroup of a conjugate of a
vertex (edge) group of the given decomposition. For our purposes
this result is particularly useful when we have two different
diagrams for $W$ so that visual decompositions with respect to
the two diagrams can be played against one another.

\medskip

\noindent {\bf (3) Coxeter Quotients.} Suppose $(W,S)$ is a
Coxeter system with diagram $V$. If $T\subset W$ then let $N(T)$
be the normal closure of $T$ in $W$. If $T\subset S$ then
$W/N(T)$ is a Coxeter group with diagram obtained from $V$ by
removing the vertices of $T$ and all vertices that connect to a
vertex of $T$ by an path with all odd labeled edges. In this
paper, we often consider a diagram for an even Coxeter system
$(W,S)$ and another diagram $V'$ for the system $(W,S')$ where
$V'$ may have odd labeled edges. Our Proposition \ref{P7}
describes a 1-1 correspondence between the set of edges with
label $>2$  in $V$ and those edges with label $>2$ in $V'$. If an
edge $[xy]$ of $V'$ has odd label and $[xy]$ corresponds to the
edge $[ab]$ of $V$, then in fact, the cyclic group $\langle
xy\rangle$ is conjugate to the group $\langle (ab)^2\rangle$. A
diagram for $W/N(xy)$ is obtained from $V'$ by collapsing the
edge $[xy]$ to a point. If $[xyu]$ is triangle then our Lemma
\ref{Tri} states that $[xu]$ and $[yu]$ are labeled 2. Hence
additional ``collapsing" in $V'$ is not generated by the
collapsing of $[xy]$ (see Lemma \ref{Quo}).  A diagram for
$W/N(xy)$ is obtained from $V$ by changing the label on $[ab]$ to
2. In this way, we can be sure that $W/N(xy)$ is an even Coxeter
group with a diagram that preserves potentially desirable aspects
of $V'$. Other quotients of diagrams for $(W,S)$ and $(W,S')$ are
obtained when we find subsets $\sigma \subset S$ and
$\sigma'\subset S'$ such that $\langle \sigma\rangle$ and
$\langle \sigma'\rangle$ are conjugate. If $f:\langle
\sigma'\rangle \to \mathbb Z_2 (=\{-1,1\})$ is a homomorphism,
and $N$ is the normal closure in $W$ of $ker(f)$, then Lemma
\ref{L29} describes how to obtain an even diagram for $W/N$ from
$V$. Understanding how quotients of $W$ correspond to quotients
of two different diagrams for $W$ is crucial to the success of
our arguments in this paper.

\section{A Reduction and Outline}

The main theorem can be reduced to a combinatorial fact.

\begin{theorem} \label{T2}
Suppose $W$ is a finitely generated even Coxeter group and $V'$
is a diagram for $W$. If $l$ is a simple (does not cross itself)
edge loop in $V'$ of length $\geq 4$ and containing an odd
labeled edge, then there is an edge of $V'$ containing two
non-consecutive vertices of $l$. (I.e. there is a shortcut in
$l$.)
\end{theorem}

In Section 6, we prove Theorem \ref{T2} for loops of length 4 and
in Section 7 we show all other cases reduce to the length 4 case,
finishing the proof of Theorems \ref{T1} and \ref{T2}. In this
section we prove:

\begin{proposition} \label{P3}
Suppose $W$ is a finitely generated even Coxeter group and $V'$ is
a diagram for $W$ with odd labeled edge $[xy]$ such that every
simple closed edge loop containing $[xy]$ of length $\geq 4$ has
a shortcut, then after a twist, a triangle containing $[xy]$ can
be replaced by an edge with even label. (In particular, the main
theorem can be reduced to Theorem \ref{T2}.)
\end{proposition}

\noindent {\bf Proof:} Our hypothesis and Lemma \ref{Tri} (below)
imply:

\begin{lemma}\label{L4}
Every path in $V'$ from $x$ to $y$ either contains the edge $[xy]$
or intersects $lk_2(x,y)$. $\bullet$
\end{lemma}

Suppose $U$ is the union of all components $C$ of $V'-(\{x,y\}\cup
lk_2(x,y))$ such that there is an edge from $x$ to $C$. Then by
Lemma \ref{L4}, there is no edge from $y$ to $U$. If there is a
vertex $t\in V'-(\{x,y\}\cup U)$ that connects to $U$ by an edge,
then $t\in lk_2(x,y)$ and so $U$ can be twisted around $\{x,y\}$,
to form the diagram $\hat V'$ for $W$. Note that after twisting,
vertices of $U$ that were connected to $x$ are replaced in $\hat
V'$  by vertices that connect to $y$ instead. I.e. in $\hat V'$,
each edge (other than $[xy]$) containing $x$ has its other vertex
in $lk_2(x,y)$.

\begin{proposition} \label{P5}
Suppose $V$ is an even diagram for the Coxeter group $W$, $[xy]$
is an odd edge in a diagram $V'$ for $W$ and every edge $[xc]$ of
$V'$ for $c\not =y $ is such that $c\in lk_2(x,y)$. Then there
exists a vertex $u \in lk_2(x,y)$ such that if $[uc]$ is an edge
with $c\not\in \{x,y\}$, then $[uc]$ has label 2, $c\in
lk_2(x,y)$, and any simplex $\sigma '$ of $V'$ that contains $u$
(respectively $x$ and $y$) and is such that $\langle
\sigma'\rangle$ is conjugate to $\langle \sigma\rangle$ for
$\sigma$ a simplex of $V$, contains $\{x,y\}$ (respectively $u$).
\end{proposition}

This proposition implies that the triangle $[xyu]$ of $V'$ can be
replaced by the edge $[yb]$ with label $2$ times the order of
$xy$, finishing Proposition \ref{P3}. More specifically, form the
diagram $\hat V'$ from $V'$ by removing the vertices $x$ and $u$,
adding a vertex $b$ and edge from $y$ to $b$ with label 2 times
the order of $xy$, and for each vertex $c$ of $lk(x)-
\{u,y\}=lk(u)-\{x,y\}$, add an edge labeled 2 from $c$ to $b$.
Then $\hat V'$ is a diagram for $W$ satisfying the conclusion of
Proposition \ref{P3}. $\bullet$

\medskip

The proof of Proposition \ref{P5} requires the development of some
basic results and is postponed until Section 5.

\section{Classifying the Rigid Even Coxeter Groups}

In this section we develop several important tool lemmas and
prove Theorem \ref{T26}, a classification of the rigid even
Coxeter groups. Through the remainder of the paper we rely on
\cite{Bourbaki} as a reference for basic facts about Coxeter
groups.

If $V$ is a diagram for a Coxeter group $W$, then a simplex
$\sigma$ is {\it spherical} if $\langle \sigma \rangle$ is a
finite subgroup of $W$ and $\sigma$ is {\it maximal spherical} if
$\sigma$ is spherical and properly contained in no other
spherical simplex. Maximal spherical simplices of $V$ give (up to
conjugation), the maximal finite subgroups of $W$. Hence if $V$
and $V'$ are diagrams for $W$ and $\sigma$ is a maximal spherical
simplex of $V$, then there is a maximal spherical simplex $\sigma
'$ of $V'$ such that $\langle \sigma\rangle$ is conjugate to
$\langle \sigma'\rangle$.

For any Coxeter diagram $V$, let $T(V)$ be the product of all
edge labels of $V$.

\medskip

\noindent {\bf Remark 1} For any integer $k$, $D_{2(2k+1)}\equiv
\langle u,v : u^2,v^2,(uv)^{2(2k+1)}\rangle =\langle u,
vuv\rangle \times \langle (uv)^{2k+1}\rangle \equiv
D_{2k+1}\times \mathbb Z_2$. If $n$ is not of the form $2(2k+1)$
then $D_n$ is irreducible.

\begin{lemma} \label{L6}Suppose the group $G$ decomposes as direct
products $\Pi_{i=1}^qA_i=\Pi _{i=1}^q B_i$ where each $A_i$ and
$B_i$ is either $\mathbb Z_2$ or $D_k$ for $k\not=2(2m+1)$, (i.e.
$D_k$ is an irreducible dihedral group). If $B_i=\langle
x,y:x^2,y^2,(xy)^n\rangle$ then there exists a unique integer $j$
such that $A_j=\langle u,v;u^2,v^2,(uv)^n\rangle$ and

\item{(i)} For odd $n$, $xy=(uv)^p$ and $\langle xy\rangle
=\langle uv\rangle$

\item{(ii)} For even $n$, $xy=(uv)^pt$ where $t$ is order 2 and
commutes with $u$ and $v$ and $\langle (xy)^2\rangle =\langle
(uv)^2\rangle$.
\end{lemma}

\noindent {\bf Proof:} In either case, for all $t\in G$,
$txyt^{-1}=(xy)^{\pm 1}$. Say $xy=a_1\cdots a_q$ where $a_i\in
A_i$. Since $xy$ does not have order 2, we may assume that $a_1$
has order greater than 2. Say $A_1=\langle u,v:
u^2,v^2,(uv)^m\rangle$. Then $a_1=(uv)^p$ and $(xy)^{\pm
1}=u(xy)u=a_1^{-1}a_2\cdots a_q$. As $a_1\not =a_1^{-1}$, we must
have $a_i=a_i^{-1}$ for all $i\geq 2$. If $xy$ has odd order, then
each $a_i$ is trivial or has odd order. In this case, $a_i=1$ for
all $i\geq 2$ and $xy=a_1$. If $n$ is even then $xy=a_1t$ where
$t$ has order 2 and commutes with $u$ and $v$.

In any case, the cyclic group $\langle xy\rangle$ is normal in
$G$ and the quotient of $G$ by $\langle xy\rangle$ has
irreducible decomposition obtained from $\Pi _{j=1}^qB_j$ by
replacing $B_i$ by $\mathbb Z_2$. Suppose $n$ is odd, and
$a_1=(uv)^p$. By the Krull-Schmidt theorem (see \cite{Ro}),
$\langle (uv)^p\rangle$ must have index 2 in $A_1$ and so $\langle
(uv)^p\rangle=\langle uv\rangle$ as desired. This implies that
$uv$ and $xy=(uv)^p$ have the same order and so $m=n$.

Now suppose $n$ is even. In this case, $(xy)^2=a_1^2=(uv)^{2p}$.
The quotient of $G$ by the normal subgroup $\langle (xy)^2\rangle
$ has decomposition obtained from $\Pi _{j=1}^qB_j$ by replacing
$B_i$ by $\mathbb Z_2\times \mathbb Z_2$. Hence the quotient of
$A_1$ by $\langle (uv)^{2p}\rangle$ is $\mathbb Z_2\times \mathbb
Z_2$. If $k$ is odd, then $D_k$ does not map onto $\mathbb
Z_2\times \mathbb Z_2$ so $uv$ has even order. If $(uv)^2\not \in
\langle (uv)^{2p}\rangle$, then 1, $u$, $uv$, $(uv)^2$ and
$(uv)^3$ would represent 5 different $(uv)^{2p}$-cosets of
$\langle u,v\rangle$ which is impossible.  Hence $\langle
(uv)^{2p}\rangle=\langle (uv)^2\rangle$. This implies that $xy$
and $uv$ have the same order and so $m=n$.

In either case $\langle (xy)^2\rangle =\langle (uv)^2\rangle $.
The uniqueness follows by construction. $\bullet$

\medskip

If a triangle $[xyz]$ in a diagram $V$ for a Coxeter group has
edge labels $(a,b,c)$, then $\langle x,y,z\rangle$ is finite if
and only if ${\frac 1 a} +{\frac 1 b} +{\frac 1 c} >1$. A result
of Tits, implies that if $A$ is a finite subgroup of a Coxeter
group $W$, and $(W,S)$ is a Coxeter system, then there is a
subset $T$ of $S$ such that $\langle T\rangle$ is finite and $A$
is a subgroup of a conjugate of $\langle T\rangle$. If $(W,S)$ is
even, and $T\subset S$ is such that $\langle T\rangle$ is finite,
then $\langle T\rangle$ decomposes as a direct product of groups
each factor of which is dihedral or $\mathbb Z_2$. It is
straightforward to see that the finite triangle groups (2,3,3),
(2,3,4) and (2,3,5) are not isomorphic to a subgroup of a
dihedral group $D_n$. The Krull-Schmidt theorem then implies that
these finite triangle groups are not subgroups of an even Coxeter
group.

If $V$ is a diagram for a Coxeter system $(W,S)$ and $[ab]$ is an
even labeled edge of $V$, then a simple examination of
presentations shows that $W/N((ab^2))$ is a Coxeter group with
diagram obtained from $V$ by changing the label of $[ab]$ to 2 .
The subgroup $\langle S-\{t\}\rangle$ for $t\in \{a,b\}$ injects
under the quotient.

\begin{proposition} \label{P7}
Suppose $W$ is a finitely generated even Coxeter group with
diagrams $V$ and $V'$ (not necessarily even).

There is a unique bijection $\alpha$ between the edges $[xy]$ of
$V'$ with label $>2$ and the edges $[ab]$ of $V$ with label $>2$
such that if $\alpha ([xy])=[ab]$, then exactly one of the
following holds.
\item{(i)} $[xy]$ and $[ab]$ are labeled $2k+1$, $xy$ is conjugate
to $(ab)^p$ for some integer $p$  and the cyclic group $\langle
xy\rangle$ is conjugate to $\langle ab\rangle$
\item{(ii)} $[xy]$ is labeled $2k+1$
and $[ab]$ is labeled $2(2k+1)$, $xy$ is conjugate to $(ab)^{2p}$
and the cyclic group $\langle xy\rangle$ is conjugate to $\langle
(ab)^2\rangle$.
\item{($ii'$)} $[xy]$ is labeled $2(2k+1)$
and $[ab]$ is labeled $2k+1$, $ab$ is conjugate to $(xy)^{2p}$ and
the cyclic group $\langle (xy)^2\rangle$ is conjugate to $\langle
ab\rangle$.
\item{(iii)} $[xy]$ and $[ab]$ are labeled $2(2k+1)$, $(xy)^2$ is
conjugate to $(ab)^{2p}$ and the cyclic group $\langle
(xy)^2\rangle$ is conjugate to $\langle (ab)^2\rangle$.
\item{(iv)} $[xy]$ and $[ab]$ are labeled $4n$, $(xy)^2$ is conjugate to
$(ab)^{2p}$ and the cyclic group $\langle (xy)^2\rangle$ is
conjugate to $\langle (ab)^2\rangle$.
\end{proposition}

\noindent {\bf Proof:} Observe that the Proposition can be
reduced to showing parts $(ii)$, $(iii)$ and $(iv)$ are valid
when all edges of $V$ have even labels. So we make that
assumption. Let $\sigma ' $ be a maximal spherical simplex of $V'$
containing $[xy]$ and $\sigma$ the maximal spherical simplex of
$V$ such that $\langle \sigma '\rangle$ is conjugate to $\langle
\sigma \rangle$. By conjugating, we assume $\langle \sigma
'\rangle = \langle \sigma \rangle$. Let $e'_1,\ldots ,e'_n$ be
the edges of $\langle \sigma '\rangle$ not labeled by 2. Then
$\langle \sigma '\rangle$ naturally decomposes as a direct
product $A'\cong \Pi _{i=1}^m A_i'$ and for $i\in \{1,\ldots ,
n\}$, $A_i'$ is the dihedral group $D_{k_i}$ where $k_i$ is the
label of $e_i'$ when $e_i'$ has odd label or label a multiple of
$4$, and $k_i$ is half the label of $e_i'$ if $e_i'$ has label
two times an odd. So if $e_i'=[xy]$ has label $2$ times an odd
integer $q$, then $q=k_i$ and $D_{2q}=\langle x,y\rangle$
decomposes as in Remark 1 as $\mathbb Z_2\times D_q = \langle
(xy)^q \rangle \times \langle x,yxy\rangle $. All other $A_j'$ are
copies of $\mathbb Z_2$. Similarly decompose $\langle \sigma
\rangle$ as the direct product $A\cong \Pi_{i=1}^m A_i$. (Note
that the number of factors in the decompositions of $\langle
\sigma \rangle$ and $\langle \sigma '\rangle$ are the same and
there is a bijection $\phi$ of the set of $A_i'$ to the set of
$A_i$ such that $\phi (A_i')$ is isomorphic to $A_i'$, by
Kurll-Schmidt.)

Now apply Lemma \ref{L6} to get a map $\alpha$ from the edges with
label $>2$ of $V'$ to those of $V$. If $[ab]$ and $[cd]$ are
distinct edges of $V$ with labels $>2$,  then $\langle
(ab)^2\rangle$ is not conjugate to $\langle (cd)^2\rangle$ since
$\langle (cd)^2\rangle$ injects under the quotient of $W$  by
$N((ab)^2)$ (the normal closure of $(ab)^2$ in $W$). Hence there
is exactly one choice for $\alpha$. By considering maximal
spherical simplices in $V$, we see that $\alpha$ is onto.

It remains to show that $\alpha$ is injective. First we prove the
several lemmas.

\begin{lemma} \label{L8}
Suppose $V$ and $V'$ are diagrams for an even Coxeter group $W$
and $V$ is even. Then no triangle of $V'$, containing an edge with
odd label, contains two edges corresponding to distinct edges of
$V$. Equivalently, if $[xyz]$ is a triangle of $V'$ and some edge
of $[xyz]$ has odd label and corresponds to the edge $[ab]$ of
$V$ then either all edges of $[xyz]$ have odd label and
correspond to $[ab]$, or two edges of $[xyz]$ are labeled 2, or
two edges of $[xyz]$ have odd labels and correspond to $[ab]$ and
the other edge is labeled 2.
\end{lemma}

\noindent {\bf Proof:} Suppose $[xy]$ is labeled $2k+1$ and
corresponds to the edge $[ab]$ of $V$. Also assume that $[xz]$ is
labeled $m>2$ and corresponds to $[cd]\not =[ab]$. Then either
$[zy]$ is labeled 2, or $[zy]$ corresponds to $[ab]$, or $[zy]$
corresponds to $[cd]$, or $[zy]$ corresponds to $[ef]\not \in
\{[cd], [ab]\}$. Observe that $N(xy)=N((ab)^2)$. We consider the
quotient map $q:W\to W/N(xy)$. Now, $\langle c,d\rangle $ injects
under $q$, but if $[zy]$ is labeled 2, then $(xz)^2\in ker(q)$.
Either $\langle xz\rangle$ is conjugate to $\langle
(cd)^2\rangle$  and $(cd)^2$ has odd order, or $\langle
(xz)^2\rangle$ is conjugate to $\langle (cd)^2\rangle $. In the
first case, $(cd)^4$ is in $ker(q)$ which is impossible. In the
second case, $(cd)^2$ is in $ker (q)$ which is also impossible.

If $[zy]$ corresponds to $[ab]$ then $[zy]$ has label $2k+1$ or
$2(2k+1)$. This label cannot be $2(2k+1)$ since otherwise,
$N((zy)^2)=N((ab)^2)=N(xy)$, but $xy\not \in N((zy)^2)$. This
label cannot be $2k+1$, since otherwise, $N(zy)=N((ab)^2)=N(xy)$.
But, $xy\in N(yz)$ implies $xz\in N(yz)$. This implies $(cd)^2\in
N(yz)=N((ab)^2)$. But $(cd)^2\not \in N((ab)^2)$.

If $[yz]$ corresponds to $[cd]$, then assume that $[cd]$ has
label $2n$. The labels of $[xz]$ and $[yz]$ are $2n$ or $n$. If
$[xz]$ is labeled $2n$, then $N((xz)^2)=N((cd)^2)$ is equal to
$N((yz)^2)$ or $N(yz)$, but $(yz)$, $(yz)^2\not \in N((xz)^2)$.
Hence $[xz]$ and similarly $[yz]$ is labeled $n$.  Then
$N(xz)=N((cd)^2)=N(zy)$. If $zy\in N(xz)$, then $xy\in N(xz)$
implying $N(xy)\equiv N((ab)^2)<N(xz)=N((cd)^2)$. But,
$(ab)^2\not \in N((cd)^2)$.

If $[yz]$ corresponds to $[ef]$, then let $q:W\to
W/N(xy)=W/N((ab)^2)$ be the quotient map.  We have $q(zx)=q(zy)$.
If the labels of $[xz]$ and $[zy]$ are odd, then $q(\langle
(cd)^2\rangle )=q(\langle xz\rangle )=q(\langle yz\rangle
)=q(\langle (ef)^2\rangle )$. But this is impossible as $q(\langle
(cd)^2\rangle )\not =q(\langle (ef)^2\rangle )$ in $W/N((ab)^2)$.
If the labels of $[xz]$ and $[zy]$ are even, then $q(\langle
(cd)^2\rangle )=q(\langle (xz)^2\rangle )=q(\langle (zy)^2\rangle
)=q(\langle (ef)^2\rangle )$ which is again impossible. If the
label of $[xz]$ is odd and the label of $[zy]$ is even then
$q(\langle (cd)^4\rangle )=q(\langle (xz)^2\rangle )=q(\langle
(zy)^2\rangle )=q(\langle (ef)^2\rangle )$. Again, $q(\langle
(cd)^4\rangle )\not =q(\langle (ef)^2\rangle )$ in $W/N((ab)^2)$.
Similarly if the label of $[xz]$ is even and the label of $[zy]$
is odd.

We conclude that each edge of $[xyz]$ with label $>2$ corresponds
to $[ab]$. If $[xz]$ has even label $>2$, then
$N((xz)^2)=N(((ab)^2)=N(xy)$ but $xy\not \in N((xz)^2)$. $\bullet$

\medskip

Now we improve Lemma \ref{L8}

\begin{lemma} \label{Tri}
Suppose $W$ is a finitely generated even Coxeter group. If $V'$
is a diagram for $W$ with odd labeled edge $[xy]$ then any
triangle of $V'$ containing $[xy]$ has two edges labeled 2.
\end{lemma}

\noindent {\bf Proof:} Suppose $V'$ is a diagram for a minimal
(with respect to $T(V')\equiv$ the product of all edge labels of
$V'$) counterexample to the Lemma and $V$ is an even diagram for
$W$. By Lemma \ref{L8}, we may assume $[xyz]$ is a triangle of
$V'$ with at least 2 odd labeled edges corresponding to $[ab]$ in
$V$. Let $\sigma '$ be a maximal simplex containing $[xyz]$. Then
by \cite{MT}, $\langle \sigma '\rangle$ is conjugate to $\langle
\sigma \rangle$ for $\sigma $ a maximal simplex of $V$. By
minimality, $V'=\sigma'$, so we have $V'$ and $V$ are complete.
Let $\delta'$ be the set of all vertices in $V'$ that belong to
an edge that corresponds to $[ab]$. If $[st]$ is an edge of $V '$
such that $s\in \delta'$ and $t\not \in \delta'$, then by Lemma
\ref{L8}, $[st]$ has label 2. I.e. $V'$ decomposes as $\langle
\delta '\rangle\times \langle V'-\delta'\rangle$. If any edge
$[st]$ of $V'$ were labeled by an even $>2$, then $W/N((st)^2)$
would be a smaller counterexample. If $s,t\in V'-\delta'$ and
$[st]$ is an edge with odd label, then $N(st)\subset
N(V'-\delta')$ and so $\langle \delta'\rangle$ injects under the
quotient map $W\to W/N(st)$. But then $W/N(st)$ is a smaller
counterexample. Hence $[ab]$ is the only edge of $V$ not labeled
2. I.e. $W$ is finite. This is impossible as a $(2,n,n)$ triangle
group is infinite unless $n=2$ or $3$. We previously ruled out
$(2,3,3)$ triangle groups as subgroups of even Coxeter groups.
$\bullet$

\medskip

As a direct consequence of the previous lemma we have the
following.

\begin{lemma}\label{Quo}
Suppose $W$ is a finitely generated even Coxeter group and $V'$
is a diagram with odd labeled edge $[xy]$. Then the diagram for
$W/N(xy)$ obtained from $V'$ by collapsing the edge $[xy]$ is
such that no other edge of $V'$ is collapsed and the only edges
of $V'$ that are identified are those in a triangle containing
$[xy]$. $\bullet$
\end{lemma}

\medskip

It is now elementary to show the injectivity of $\alpha$. If
distinct edges $[xy]$ and $[zw]$ of $V'$ correspond to $[ab]$ in
$V$, then a quotient argument easily implies neither $[xy]$ nor
$[zw]$ has even label. But then $N(xy)=N((ab)^2)=N(zw)$. But by
Lemma \ref{Quo}, $zw\not\in N(xy)$. $\bullet$

\medskip

\noindent {\bf Example 1.} The element $cba$ of the group
$\langle a,b,c: a^2=b^2=c^2=(ab)^3=(bc)^3=(ac)^2=1\rangle$
conjugates $b$ to $c$ and $a$ to $b$. Indicating the need for a
more sophisticated version of Proposition \ref{P7} in a more
general setting.

\medskip

The Deletion Condition for Coxeter groups implies the following
Lemma.

\begin{lemma} \label{L11}
Suppose $(W,S)$ is a Coxeter system, $\Gamma$ the Cayley graph of
$W$ with respect to $S$ and $T\subset S$. If $u$ and $v$ are
vertices of $\Gamma$ (i.e. elements of $W$), then there is a
unique closest vertex $w$ of the coset $v\langle T\rangle$ to
$u$. Furthermore, if $\alpha$ is a geodesic from $u$ to $w$, and
$\beta$ is a geodesic at $w$ in the letters of $T$, then $\alpha
\beta$ is geodesic. $\bullet$
\end{lemma}

\begin{proposition} \label{P12}
Suppose $(W,S)$ is an even Coxeter system, $a,b\in S$ and $ab$
has finite order $>2$. If $y\in W$ is such that $y$ conjugates
$(ab)^2$ to $(ab)^{\pm 2}$ then $y$ can be written geodesically
as $uv$ where $u\in \langle a,b\rangle$ and $v\in lk_2(a,b)$.
\end{proposition}

\noindent {\bf Proof:} We first show that $y$ conjugates $\langle
a,b\rangle$ to itself. We have $(ab)^2\in \langle a,b\rangle\cap
y\langle \langle a,b\rangle y^{-1}=v\langle T\rangle v^{-1}$ for
$T\subset \{a,b\}$ and $v\in \langle a,b\rangle$. If $T$ is a
single element, then $(ab)^2 $ is conjugate to $a$ or $b$. This is
impossible as $(ab)^2$ has even length. Hence $T=\langle
a,b\rangle$ and so $\langle a,b\rangle=y\langle a,b\rangle
y^{-1}$.

Let $\Gamma$ be the Cayley graph of $W$ with respect to $S$.
Write $y=x_1y_1x_2$ where $x_i\in \langle a,b\rangle$ and $y_1 $
is the shortest element of the double coset $\langle a,b\rangle
y\langle a,b\rangle$. We show that $y_1$ commutes $a$ and $b$. If
$\alpha$ is a geodesic in $\Gamma$ from $1$ to $y_1$ then
$a\alpha$, $b\alpha$, $\alpha a$ and $\alpha b$ are geodesic by
the choice of $y_1$. Hence by Lemma \ref{L11} if $\beta _1$ and
$\beta _2$ are geodesic paths at 1 and $y_1$ respectively, in the
letters $a,b$, then the paths $(\beta _1^{-1}, \alpha )$ and
$(\alpha , \beta _2)$ are geodesic. Now, since $y_1ay_1^{-1}$ and
$y_1by_1^{-1}$ are in $\langle a,b\rangle$ they must both be of
length 1. I.e. (since $(W,S)$ is even) $y_1$ commutes with $a$
and $b$. Furthermore, $y=x_1x_2y_1$. Results in \cite{BH} and
\cite{Bo} implies $y_1$ is a product of an element of $\langle
a,b\rangle$ and an element of $lk_2(a,b)$.  $\bullet$

\medskip

\noindent {\bf Remark 2} Each part of Proposition \ref{P7}
concludes that $xy$ or $(xy)^2$ is conjugate to $(ab)^p$ or
$(ab)^{2p}$. Hence if $c$ is the conjugating element, then
$cabc^{-1}$ commutes with $(xy)^2$ and by Proposition \ref{P12},
$cabc^{-1}=uv$ for $u\in \langle x,y\rangle$ and $v\in
lk_2(x,y)$. Similarly, $c^{-1}xyc$ is conjugate to $u'v'$ for
$u\in \langle a,b\rangle$ and $v'\in lk_2(a,b)$. Hence, parts
$(ii')$, $(iii)$ and $(iv)$ of Proposition \ref{P7} can be
improved to say: $(xy)$ is conjugate to $(ab)^pt$, where $t$
commutes with $a$ and $b$ and $t^2=1$. Parts $(ii)$, $(iii)$ and
$(iv)$ can be improved to say: $(ab)$ is conjugate to $(xy)^qs$,
where $s$ commutes with $x$ and $y$ and $s^2=1$.

\begin{proposition} \label{P13}
Suppose $(W,S)$ is an even Coxeter system with diagram $V$,
$(W,S')$ is another Coxeter system with diagram $V'$, $[xy]$ and
$[yz]$ are distinct edges of $V'$ with labels $>2$. If $[ab]$ and
$[cd]$ are edges of $V$ which correspond to $[xy]$ and $[yz]$
respectively, then $\{a,b\}\cap \{c,d\}$ contains exactly one
element. I.e. The edges $[ab]$ and $[cd]$ share exactly one
vertex.
\end{proposition}

\noindent {\bf Proof:} Suppose that $[xy]$ has label $k>2$ and
$[yz]$ has label $m>2$. Observe that $[ab]$ is labeled $l$ where
$l=k$  or $l=2k$ and $[cd]$ by $n$ where $n=m$ or $n=2m$. The set
$\{a,b\}\cap \{c,d\}\not =\{a,b\}$ by the uniqueness of pairing
of Proposition \ref{P7}. Assume that $\{a,b\}\cap
\{c,d\}=\emptyset$. Now $xy$ or $(xy)^2$ is equal to
$w(ab)^{2p}w^{-1}$ for some $w\in W$. By conjugation we may
assume that either $yz$ or $(yz)^2$ is equal to $(cd)^{2r}$.

\begin{lemma}\label{L14}
Not both $a$ and $b$ commute with both $c$ and $d$.
\end{lemma}

\noindent {\bf Proof:} Otherwise there is a maximal spherical
simplex $\sigma $ (i.e. $\langle \sigma \rangle$ is finite)
containing $a,b,c$ and $d$. The group $\langle \sigma \rangle$ is
conjugate to $\langle \sigma '\rangle$ for $\sigma '$ a simplex
of $V'$ containing $x,y$ and $z$, which is impossible as $\langle
x,y,z\rangle$ is not finite. $\bullet$

\medskip

Consider the retraction $\alpha :W\to \langle
a,b,c,d\rangle\equiv G$ with kernel $N(S-\{a,b,c,d\})$.

\begin{lemma} \label{L15}
Either $a$ or $b$ is an element of $lk_2(c,d)$ and either $c$ or
$d$ is an element of $lk_2(a,b)$.
\end{lemma}

\noindent {\bf Proof:} Let $\alpha (x)\equiv \bar x$, $\alpha
(y)\equiv \bar y$, $\alpha (z)\equiv \bar z$ and $\alpha
(w)\equiv \bar w$. As $(ab)^{2p}=w^{-1}(xy)w$ or $w^{-1}(xy)^2w$,
Proposition \ref{P12} implies $\bar w^{-1}\bar y\bar w=u_1v_1$,
where $u_1\in \langle a,b\rangle$ and $v_1\in lk_2(a,b)$ (where
$lk_2$ is taken in $\langle a,b,c,d\rangle$), and $\bar y=u_2v_2$
where $u_2\in \langle c,d\rangle $ and $v_2\in lk_2(c,d)$. Note
that $v_1=c$ or $d$ or $1$ and $v_2=a$ or $b$ or $1$. If $v_2=1$,
then $\bar y\in \langle c,d\rangle$. The element $\bar w^{-1}\bar
y\bar w$ conjugates $(ab)^{2}$ to $(ba)^{2}$. As $\bar y$ ($\in
\langle c,d\rangle$) is in the kernel of the retraction of
$\langle a,b,c,d\rangle$ to $\langle a,b\rangle $ (with kernel
$N(\{c,d\})$), this is impossible unless $(ab)^{2} =(ba)^{2}$. If
$ab$ has order 4, then $xy$ has order 4 and Remark 2 implies
$xyt=w(ab)^qw^{-1}$ where $t$ has order 2 and commutes with $x$
and $y$. In this case we see that $w^{-1}\bar yw$ conjugates $ab$
to $ba$. Again this is impossible as $\bar y\in \langle
c,d\rangle$ is in the kernel of a retraction of $\langle
a,b,c,d\rangle$ to $\langle a,b\rangle$ and $ab\not =ba$.
Similarly, $v_1\not =1$. $\bullet$

\medskip

Without loss, we assume that $a\in lk_2(c,d)$ and $c\in
lk_2(a,b)$. (See Figure 1.)

\medskip

\centerline{\bf Figure 1}

\medskip

Recall, $T(V')$ is the product of the edge labels of $V'$. Assume
that $V'$ is a minimal (under $T$) counterexample to Proposition
\ref{P13}.

\begin{lemma}\label{L16}
There is no edge $[bd]$.
\end{lemma}

\noindent {\bf Proof:} Otherwise, $[bd]$ has label $>2$ by Lemma
\ref{L14}. If the edge $[st]$ of $V'$ corresponds to $[bd]$ and
$st$ has even label $2k$, then $[bd]$ has label $2k$ and
$W/N((bd)^2)$ is a smaller counterexample to our proposition.
Hence we may assume that $[st]$ has odd label. By Lemma
\ref{Tri}, $\{s,t\}\not =\{x,z\}$, and by Lemma \ref{Quo},
$W/N(st)=W/N((bd)^2)$ is a smaller counterexample. $\bullet$

\begin{lemma} \label{L17}
There is no edge $[xz]$.
\end{lemma}

\noindent {\bf Proof:} Otherwise, let $\sigma '$ be a maximal
simplex of $V'$ containing the triangle $[xyz]$. By \cite{MT},
there is a simplex $\sigma$ of $V$ such that $\langle \sigma
'\rangle$ is conjugate to $\langle \sigma \rangle$. Thus,
$\{a,b,c,d\}\subset \sigma$ contradicting Lemma \ref{L16}.
$\bullet$

\begin{lemma} \label{L18}
In $V'$, the only edges not labeled 2 are $[xy]$, $[yz]$.
\end{lemma}

\noindent {\bf Proof:} By the minimality assumption, there are no
edges of $V'$ with even label $>2$. If $[st]$ is distinct from
$[xy]$ and $[yz]$, and with odd label, then the quotient of $W$
by $N(st)$ gives a smaller counterexample by Lemma \ref{Quo}.
$\bullet$

\medskip

Let $S'$ be the vertex set of $V'$. Let $\lambda$ be the
retraction of $W$ to $\langle x,y,z\rangle$ with kernel
$N(S'-\{x,y,z\})$.  Observe that the groups $\langle b,c\rangle$,
$\langle a,d\rangle$, and $\langle a,c\rangle$ inject under
$\lambda$, for otherwise $(ab)^2$ or $(cd)^2$ is in $ker
(\lambda)$ (see Figure 1) which they are not.

Now $\langle x,y,z\rangle =\langle x,y\rangle \ast _{\langle
y\rangle}\langle y,z\rangle$.  A maximal simplex $\sigma$ of $V$
containing $\{a,b,c\}$ is such that $\langle \sigma \rangle$ is
conjugate to a $\langle \sigma'\rangle$ for $\sigma '$ a maximal
simplex of $V'$. Hence $\lambda (\langle a,b,c\rangle )$ is a
subgroup of a conjugate of $\langle x,y\rangle$ or $\langle
y,z\rangle$. As $\lambda ((ab)^2)$ is in the kernel of the
quotient of $\langle x,y,z\rangle $ by $N(xy)$, and $\langle
y,z\rangle$ injects under this quotient, $\lambda (\langle
a,b,c\rangle )$ is a subgroup of a conjugate of $\langle
x,y\rangle$. Similarly, $\lambda (\langle a,c,d\rangle )$is a
subgroup of a conjugate of $\langle y,z\rangle$. Hence (simply
consider edge and vertex stabilizers of the Bass-Serre tree for
$\langle x,y\rangle \ast _{\langle y\rangle}\langle y,z\rangle$)
$\lambda (\langle a,c\rangle )$ is a subgroup of a conjugate of
$\langle y\rangle$. But this is impossible as $\langle y\rangle$
has order 2. $\bullet$

\medskip

\noindent {\bf Example 2} Three diagrams for an even Coxeter
group are shown in Figure 2. Isomorphisms between presentations
determined by these diagrams are given by:

\noindent $a\to a$, $c\to c$, $d\to d$, $x\to bab$ and
$y\to(ab)^3$ (a triangle/edge exchange) and

\noindent $a\to a$, $c\to c$, $x\to x$, $y\to y$ and $z\to
axadaxa$ (a twist of $\{d\}$ around $[ax]$).

\noindent The correspondence of Proposition 7 between edges of the
first and last diagram of Figure 2, match the adjacent edges
$[ac]$ and $[ad]$ with the non-adjacent edges $[ac]$ and $[xz]$,
respectively.

\medskip

\centerline{\bf Figure 2}

\medskip

\begin{proposition} \label{P19}
Suppose $(W,S)$ is an even Coxeter system with diagram $V$,
$(W,S')$ is another Coxeter system with diagram $V'$, $[ab]$ and
$[ac]$ are distinct edges of $V$ with labels $>2$. If $[xy]$ and
$[zv]$ are edges of $V'$ which correspond to $[ab]$ and $[ac]$
respectively, then either there is an odd edge path in $V'$
connecting a vertex of $\{x,y\}$ and a vertex of $\{z,v\}$, or
$\{x,y\}\cap \{z,v\}$ contains exactly one element.
\end{proposition}

\noindent {\bf Proof:}   Suppose there is no odd edge path in
$V'$ with one vertex in $\{x,y\}$ and one vertex in $\{z,v\}$,
$\{x,y\}\cap \{z,v\} =\emptyset$, and $V'$ is a minimal (over
$T(V')$) such counterexample to our proposition.

\begin{lemma}\label{L20}
The only edges of $V'$ not labeled 2 are $[xy]$ and $[zv]$.
\end{lemma}

\noindent {\bf Proof:} By the minimality of $T(V')$, there is no
edge of $V'$ with a label that is even and $>2$, other than
possibly $[xy]$ and $[zv]$. Also by the minimality of $T(V')$ and
Lemma \ref{Quo} there is no edge of $V'$ with odd label other than
$[xy]$ or $[zv]$. $\bullet$

\medskip

By Lemma \ref{L20} and Proposition \ref{P7}, the only edges of
$V$ not labeled 2 are $[ab]$ and $[ac]$. Let $\sigma $ be a
maximal spherical simplex of $V$ containing $\{a,b\}$ and $\sigma
'$ be a maximal spherical simplex in $V'$ such that $\bar
w\langle \sigma'\rangle \bar w^{-1}=\langle \sigma \rangle$ for
some $\bar w\in W$. Then $\{x,y\}\subset \sigma '$. If neither
$z$ nor $w$ is an element of $\sigma '$ then
$W/N(\sigma)=W/N(\sigma ')$ is a Coxeter group with diagram
obtained from $V$ (respectively $V'$) by removing the vertices of
$\sigma$ (respectively $\sigma '$). But one of these diagrams has
all edges labeled 2 and the other has an edge with label $>2$,
which is impossible. We conclude that $v$ or $z$ is an element of
$\sigma '$ and so $v$ or $z$ is an element of $lk_2(x,y)$.
Similarly, $x$ or $y$ is an element of $lk_2(z,v)$. Say $z\in
lk_2(x,y)$ and $y\in lk_2(z,v)$. See Figure 3.

\medskip

\centerline{\bf Figure 3}

\medskip

Consider the retraction $\tau$ of $W$ to $\langle a,b,c\rangle$
with kernel $N(S-\{a,b,c\})$. By Proposition \ref{P7}, $\{(xy)^2,
(zv)^2\}\cap ker (\tau)=\emptyset$. Hence
$\{x,y,z,v,zy,xz,yv\}\cap ker(\tau)=\emptyset$ (see Figure 3).
There is no edge $[bc]$ as every simplex of $V'$ (and hence every
simplex of $V$) is spherical. Observe that $\langle a,b,c\rangle =
\langle a,b\rangle\ast _{\langle a\rangle}\langle a,c\rangle$.
Note that $\tau (xy)$ has order $>2$, and is an element of a
conjugate of $\langle a,b\rangle$. Hence $\tau (xy)$ cannot be an
element of distinct conjugates of $\langle a,b\rangle$, or some
conjugate of $\langle a,c\rangle$, since otherwise, (simply
consider vertex and edge stabilizers of the Bass-Serre tree for $
\langle a,b\rangle\ast _{\langle a\rangle}\langle a,c\rangle $)
$\tau (xy)$ is an element of a conjugate of $\langle a\rangle$, an
order 2 group. Therefore, $\tau (\langle x,y,z\rangle )$ is a
subgroup of a conjugate of $\langle a,b\rangle$ and $\tau
(\langle y,z,v\rangle )$ is a subgroup of a conjugate of $\langle
a,c\rangle$. But this implies (again consider vertex and edge
stabilizers of the Bass-Serre tree for $ \langle a,b\rangle\ast
_{\langle a\rangle}\langle a,c\rangle $) that $\alpha (\langle
z,y\rangle )$ is a subgroup of a conjugate of $\langle a\rangle$.
This is impossible and the proof of the proposition is complete.
$\bullet$

\medskip

The following result is used in \cite{B}, so we cannot use
\cite{B} to simplify the proof.

\begin{lemma} \label{L21}
Suppose $V$ is an even diagram for the finitely generated Coxeter
group $W$ and $V'$ is another diagram for $W$. If $[abc]$ is a
triangle of $V$ having at least two edges with label $>2$, then
the edges of $[abc]$ with label $>2$ correspond to edges with
label $>2$ of a triangle $[xyz]$ of $V'$.
\end{lemma}

\noindent {\bf Proof:} Let $\sigma $ be a  maximal simplex of $V$
containing $[abc]$. By \cite{MT}, there is a maximal simplex
$\sigma '$ in $V'$ such that $\langle \sigma\rangle $ is
conjugate to $\langle \sigma '\rangle$. By Proposition \ref{P7}
applied to $\sigma$ and $\sigma'$ and the uniqueness conclusion
of Proposition \ref{P7} applied to $V$ and $V'$, $\sigma'$
contains the edges of $V'$ corresponding to those of $[abc]$ that
have label $>2$. First we show the edges of $\sigma '$ with label
$>2$  and corresponding to those of $[abc]$ with label $>2$ are
mutually adjacent. If not, say $[xy]$ and $[zv]$ are two such
non-adjacent edges. By Proposition \ref{P19} (applied to the even
Coxeter group $\langle \sigma \rangle$) there is an odd labeled
edge (in $\sigma '$) adjacent to $[xy]$, but this is impossible
by Lemma \ref{Tri}.

If an edge of $[abc]$ is labeled 2, we are finished. Otherwise,
the edges of $\langle \sigma '\rangle $ corresponding to those of
$[abc]$ must have even labels by Lemma \ref{Tri} and either form a
triangle or triad. If a triad is formed, then we may assume that
$[ab]$ corresponds to $[xy]$ in $V'$, $[bc]$ corresponds to
$[xz]$ and $[ac]$ corresponds to $[xv]$. Now assume that $V$ is a
minimal (with respect to $T(V)$) counterexample to the Lemma. Then
$V=\sigma$ and $V'=\sigma'$. By minimality, every edge of $V$
except $[ab]$, $[bc]$ and $[ac]$ has label 2. Similarly for $V'$.
In particular, $V$ and $V'$ are even. By conjugation, we may
assume that $\langle (xy)^2\rangle =\langle (ab)^2\rangle$. As
$x$ conjugates $(ab)^2$ to $(ba)^2$, Proposition \ref{P12} implies
that, $x=u_1t_1$ where $u_1\in \langle a,b\rangle$ and $t_1\in
lk_2(a,b)$. Similarly, $x=w_2u_2t_2w_2^{-1}$ where $w_2\in W$,
$u_2\in \langle b,c\rangle$ and $t_2\in lk_2(b,c)$ and
$x=w_3u_3t_3w_3^{-1}$ where $w_3\in W$, $u_3\in \langle
a,c\rangle$ and $t_3\in lk_2(a,c)$.

Now $x\in (\langle \{a,b\}\cup lk_2(a,b)\rangle )\cap (w_2\langle
\{b,c\}\cup lk_2(b,c)\rangle w_2^{-1})\cap (w_3\langle
\{a,c\}\cup lk_2(a,c)\rangle w_3^{-1})= v\langle T\rangle v^{-1}$
where $T\subset \{a,b\}\cup lk_2(a,b)$. Clearly, $c\not \in T$. As
no conjugate of $a$ is an element of $w_2\langle \{b,c\}\cup
lk_2(b,c)\rangle w_2^{-1}$, $a\not\in T$. Similarly, $b\not \in
T$. But then $T$ is central in $W$, implying $x$ is central, the
desired contradiction. $\bullet$

\begin{lemma} \label{L22}
Suppose $V$ is an even diagram for the finitely generated Coxeter
group $W$ and $V'$ is another diagram for $W$. If $[xyz]$ is a
triangle of $V'$ having at least two edges with label $>2$, then
the edges of $[xyz]$ with label $>2$ correspond to edges with
label $>2$ of a triangle $[abc]$ of $V$.
\end{lemma}

\noindent {\bf Proof:} By Lemma \ref{Tri} each edge of $[xyz]$
has even label. Suppose $V'$ is a minimal (with respect to
$T(V')$) counterexample. By Lemma \ref{Quo}, $V'$ contains no odd
labeled edge, and so $V'$ is even. Now apply Lemma \ref{L21}.
$\bullet$

\medskip

\noindent {\bf Remark 3} Observe in each of the last two lemmas
that if we begin with a triangle with exactly one edge labeled 2,
we do not conclude that the corresponding triangle has an edge
labeled 2. This can now be resolved by combining the last two
lemmas. More specifically, we see that if $W$ is an even Coxeter
group, $V$ and $V'$ are diagrams for $W$ and $[abc]$ is a
triangle of $V$ with only one edge with label 2, then there is a
corresponding triangle $[xyz]$ of $V'$ with only one edge labeled
2. If $[abc]$ has no edge labeled 2 then $[xyz]$ has no edge
labeled 2.

\medskip

As a direct application of Lemma \ref{L21} and Lemma \ref{L22} we
have an analogue for Lemma \ref{Tri}.

\begin{lemma}\label{L23}
Suppose $W$ is an even Coxeter group, $V$ and $V'$ are diagrams
for $W$, $[xy]$ is an odd labeled edge of $V'$ and $[ab]$ the
edge of $V$ corresponding to $[xy]$. Then any triangle containing
$[ab]$ has two edges labeled 2. $\bullet$
\end{lemma}

\begin{proposition}\label{P24}
Suppose $V$ is an even diagram for the finitely generated Coxeter
group $W$ and $V'$ is another diagram for $W$. If $[xy]$ is an
odd labeled edge of $V'$ and $[ab]$ is the corresponding edge of
$V$, then with the exception of $[ab]$, every edge adjacent to
$a$ is labeled 2 or every edge adjacent to $b$ is labeled 2.
\end{proposition}

\noindent {\bf Proof:} Suppose $[xy]$ has label $2k+1$. It
suffices to show:

\begin{lemma} \label{L25}
There is no edge path $([ca], [ab], [bd])$ in $V$ such that each
edge has label $>2$.
\end{lemma}

\noindent {\bf Proof:} Since $[xy]$ has odd label, Lemma
\ref{L23} implies the path $([ca], [ab], [bd])$ does not form a
triangle. Suppose that $[ca]$ and $[bd]$ correspond to $[uv]$ and
$[st]$ respectively in $V'$. By Proposition \ref{P19}, select a
shortest path with odd labeled edges, $e_1,\ldots ,e_n$ from
$\{u,v\}$ to $\{x,y\}$. Assume that $e_i=[x_{i-1}x_i]$ for all
$i$. If $e_j=[st]$, then a diagram $\bar V'$ for
$W/N(\{x_0x_1,\ldots , x_{j-2}x_{j-1}\})$ is obtained from $V'$
by collapsing each edge of the set $\{e_1,\ldots ,e_{j-1}\}$. The
corresponding diagram $\bar V$ for $W/N(\{x_0x_1,\ldots
,x_{j-2}x_{j-1}\})$, is obtained from $V$ by replacing each label
of an edge corresponding to one of $\{e_1,\ldots ,e_{j-1}\}$ by
2. In $\bar V'$, $[uv]$ and $[st]$ are adjacent, but $[ca]$ and
$[bd]$ are not, contradicting Proposition \ref{P13}. We conclude
that $[st]$ is not in $\{e_1,\ldots ,e_n\}$. Similarly, if
$d_1,\ldots ,d_m$ is a shortest odd labeled edge path from
$\{x,y\}$ to $\{s,t\}$, we may assume that $[uv]$ is not in
$\{d_1,\ldots ,d_m\}$. Assume that $d_i=[y_{i-1}y_i]$ for all
$i$. A diagram $\tilde V'$ is obtained for the group $\tilde
W\equiv W/N(\{x_0x_1,\ldots ,x_{n-1}x_{n}, xy, y_0y_1,\ldots
,y_{m-1}y_m\})$ by collapsing the edges $e_1,\ldots , e_n, [xy],
d_1,\ldots ,d_m$ of $V'$. Hence in $\tilde V'$, $[uv]$ and $[st]$
are adjacent. Let $\tilde V$ be the diagram for $\tilde W$
obtained from $V$ by changing the edge labels of the edges of $V$
corresponding to $e_1,\ldots ,e_n,[xy], d_1,\ldots ,d_m$ to 2. In
$\tilde V$, $[ca]$ and $[bd]$ are not adjacent, contradicting
Proposition \ref{P13}. $\bullet$ $\bullet$

\medskip

The following theorem classifies even rigid Coxeter groups.

\begin{theorem} \label{T26}
If $V$ is an even diagram for the Coxeter group $W$ then $W$ has
a diagram that is not even if and only if there is an edge $[ab]$
in $V$ with label $2(2k+1)$ for $k>0$, such that with the
exception of $[ab]$, every edge of $V$ containing $a$ is labeled
2 and if $[ac]$ is such an edge, then there is an edge $[bc]$
with label 2.
\end{theorem}

\noindent {\bf Proof:} If $[ab]$ is an edge as described in the
theorem, then $\langle a,b:a^2=b^2=(ab)^{2(2k+1)}=1\rangle$ is
isomorphic to the group $\langle
x,y,z:x^2=y^2=z^2=(xz)^2=(yz)^2=(zy)^{2k+1}=1\rangle$ by the map
extending $x\to a$, $y\to bab$ and $z\to (ab)^{2k+1}$. It is
elementary to see that the edge $[ab]$ in $V$ can be replaced by
the triangle $[xyz]$ to give a new diagram for $W$.

The proof of the converse is more delicate. Recall that $T(V)$ is
the product of all edge labels in $V$. We assume from this point
on that $V$ is a minimal (with respect to $T$) counterexample to
our theorem. Let $V'$ be a diagram for $W$ with odd labeled edge
$[xy]$. Assume that $[xy]$ corresponds to the edge $[ab]$ of $V$.

\begin{lemma} \label{L27}
With the exception of $[ab]$ every edge of $V$ is
labeled 2 (and hence $[xy]$ is the only edge of $V'$ not labeled
2.)
\end{lemma}

\noindent {\bf Proof:} If $[cd]$ is an edge labeled $n>2$ and
$\{c,d\}\cap \{a,b\}=\emptyset$, then the quotient of $W$ by
$N((cd)^2)$ is a ``smaller" counterexample. If $[ac]$ has label
$n>2$, then there is no edge $[bc]$ by Lemma \ref{L23}. Again the
quotient of $W$ by $N((ac)^2)$ is a smaller counterexample.
Similarly, there is no edge $[bc]$ with label $>2$. $\bullet$

\begin{lemma} \label{L28}
Suppose $\sigma $ and $\sigma '$ are simplices of $V$ and $V'$
respectively such that $\langle \sigma \rangle$ is conjugate to
$\langle \sigma '\rangle$. Then $\sigma $ contains $\{a,b\}$ if
and only if $\sigma '$ contains $\{x,y\}$. Also, $\sigma $
contains exactly one element of $\{a,b\}$ if and only if $\sigma
'$ contains exactly one element of $\{x,y\}$.
\end{lemma}

\noindent {\bf Proof:} The group $\langle \sigma \rangle $
(respectively $\langle \sigma '\rangle $) is non-abelian if and
only if $\{a,b\}\subset \sigma$ (respectively $\{x,y\}\subset
\sigma '$). Hence the first conclusion of the lemma follows.

The group $W/N(\sigma)$ (respectively $W/N(\sigma ')$ is abelian
iff $a$ or $b\in \sigma$ (respectively $x$ or $y\in \sigma'$).
But $W/N(\sigma )=W/N(\sigma ')$. $\bullet$

\medskip

To finish Theorem \ref{T26}, it suffices to show that $V$ cannot
contain edges $[ad]$ and $[bc]$ ($d\ne b$ and $c\ne a$)  such that
there is no edge between $b$ and $d$ and no edge between $a$ and
$c$. (There may or may not be an edge $[cd]$.) Assume otherwise.
Let $\sigma (a,d)$ and $\sigma (a,b) $ be maximal simplices of $V$
containing $\{a,d\}$ and $\{a,b\}$ respectively. Let $\sigma
=\sigma (a,d)\cap \sigma (a,b)$. Note that $a\in \sigma$, but
$\{b,c,d\}\cap \sigma=\emptyset$. By conjugation we may assume
that $\sigma'(a,b)$ is a maximal simplex of $V'$ such that
$\langle \sigma '(a,b)\rangle=\langle \sigma (a,b)\rangle$ and
that $w\langle \sigma '(a,d)\rangle w^{-1}=\langle \sigma
(a,d)\rangle$ for $\sigma '(a,d)$ a maximal simplex of $V'$ and
$w\in W$. Then $\langle \sigma \rangle =\langle \sigma
'(a,b)\rangle \cap w \langle \sigma' (a,d)\rangle w^{-1}=v\langle
T\rangle v^{-1}$, for some $v\in \langle \sigma '(a,b)\rangle$,
and $T\subset \sigma '(a,b)$. By Lemma \ref{L28} either $x$ or
$y$, but not both is an element of $T$.

Let $q$ be the retraction of $W$ to $\langle a,b,c,d\rangle$ with
kernel $N(S-\{a,b,c,d\})$. Then, $q(\langle a\rangle) =q(\langle
\sigma (a,b)\rangle \cap \langle \sigma (a,d)\rangle )=q(v\langle
T\rangle v^{-1})$. Observe that $x$ is conjugate to $y$, $\langle
xy\rangle$ is conjugate to $\langle (ab)^2\rangle$ and $q(ab)$ has
order $2(2k+1)$. Thus, $q(x)\not =1\not =q(y)$ and so $q(\langle
a\rangle )$ is conjugate to $q(\langle x\rangle)$ and $q(\langle
y\rangle)$. Hence $q(a)$ is conjugate to $q(x)$ and $q(y)$.
Similarly for $b$. This implies $q(a)$ and $q(b)$ are conjugate,
the desired contradiction. Theorem \ref{T26} is finished.
$\bullet$

\section{The Proof of Proposition 5}

\begin{lemma} \label{L29}
Suppose $\langle \sigma'\rangle=w\langle \sigma\rangle w$ for
$\sigma '$ a simplex of $V'$, $\sigma$ a simplex of $V$ and $w\in
W$. If $f:\langle \sigma '\rangle \to \mathbb Z_2\equiv \{-1,1\}$
is a homomorphism, let $N$ be the normal closure in $W$ of $ker
(f)$. Then $W/N$ is an even Coxeter group with diagram obtained
from $V$ by removing the vertices of $\sigma _1\equiv \{s\in
\sigma  :f(wsw^{-1})=1\}$ and identifying the vertices of $\sigma
-\sigma _1$.
\end{lemma}

\noindent {\bf Proof:} The kernel of $f$ is generated by $K'$,
the normal closure in $\langle \sigma '\rangle$ of $\{s\in \sigma
':f(s)=1\}\cup\{st:s,t\in \sigma ', f(s)=f(t)\not =1\}$. Hence
$N$ is the normal closure of $K'$ in $W$.

As $\langle \sigma '\rangle =w\langle \sigma \rangle w^{-1}$,
$K'$ can also be described as the normal closure in $w\langle
\sigma \rangle w^{-1}$ of $K = \{wsw^{-1}:s\in \sigma$ and $
f(wsw^{-1})=1\} \cup \{wstw^{-1}: s,t\in \sigma$ and
$f(wsw^{-1})=f(wtw^{-1})\not =1$. Now, in $W$, the normal closure
of $K'$, $K$ and $w^{-1}Kw$ are the same. $\bullet$

\begin{lemma} \label{L30}
Suppose $(W,S)$ is an even Coxeter system with diagram $V$ and
$V'$ is another diagram for $W$ with odd edge $[xy]$. There
exists a vertex $u\in V'-\{x,y\}$ such that $u$ is contained in
the intersection of all simplicies $\sigma '$ containing
$\{x,y\}$ and such that $\langle\sigma' \rangle$ is conjugate to
$\langle \sigma \rangle$ for $\sigma$ a simplex of $V$.
Furthermore if $\sigma '$ is a simplex of $V'$ containing $u$ and
such that $\langle\sigma' \rangle$ is conjugate to $\langle
\sigma \rangle$ for $\sigma$ a simplex of $V$, then
$\{x,y\}\subset \sigma'$.
\end{lemma}

\noindent {\bf Proof:} Assume $V'$ is a minimal (with respect to
$T(V')$) counterexample. Let $\delta '$ be the intersection of all
simplices $\sigma '$ of $V'$, containing $\{x,y\}$ and such that
$\langle \sigma ' \rangle$ is conjugate to $\langle \sigma
\rangle$ for some simplex $\sigma$ of $V$. As $\langle x,y\rangle$
is not an even Coxeter group, $\delta'\not =\{x,y\}$. If $u\in
\delta '-\{x,y\}$, there is no odd path in $V'$ from $u$ to $x$
or $y$, by Lemmas \ref{Tri} and \ref{Quo}. For each $u\in \delta
'-\{x,y\}$ assume there is a simplex $\beta '$ of $V'$ such that
$\langle \beta '\rangle$ is conjugate to $\langle \beta\rangle$
for $\beta$ a simplex of $V$ and $u\in \beta '$, but $\{x,y\}\not
\subset \beta'$. By intersecting, we may assume that each such
$\beta'\subset \delta '$. Select one such $\beta '$. By Lemma
\ref{L29} the map of $\langle \beta'\rangle$ to $\mathbb Z_2$ that
sends $\beta' -\{x,y\}$ to $1$ and $\beta' \cap \{x,y\}$ to $-1$
defines a smaller counterexample. $\bullet$

\medskip

Let $\delta '$ be the intersection of all simplices $\sigma '$ of
$V'$ such that $\{x,y\}\subset \sigma '$ and $\langle \sigma
'\rangle$ is conjugate to $\langle \sigma \rangle$ for some
simplex $\sigma $ of $V$. By Lemma \ref{L30}, $\delta '$ contains
a vertex $v$ such that if $v\in \sigma'$, where $\sigma '$ is a
simplex of $V'$ and $\langle \sigma '\rangle$ is conjugate to
$\langle \sigma\rangle$ for some simplex $\sigma$ of $V$, then
$\{x,y\}\subset \sigma '$. We call such a vertex $\{x,y\}$-{\it
linked} or simply {\it linked}.

It suffices to show that $\delta '$ contains a linked vertex
$v'$, such that every edge of $V'$ containing $v'$ is labeled 2.
Otherwise, assume that $V'$ is a minimal counterexample. Then
each linked vertex belongs to an edge with label $>2$. Suppose
$[st]$ is an edge of $V'$ with label $>2$ and neither $s$ nor $t$
is linked. If $[st]$ has even (respectively odd) label, then the
even Coxeter group $W/N((st)^2)$ (respectively $W/N(st)$), with
diagram $\bar V'$, obtained from $V'$ by changing the label of
$[st]$ to a 2 (respectively identifying $s$ and $t$), is a smaller
counterexample. (Note that if a vertex is not $\{x,y\}$-linked in
$V'$, then it is not $\{x,y\}$-linked in $\bar V'$, and an
$\{x,y\}$-linked vertex of $V'$ may not be $\{x,y\}$-linked in
$V'$.) Hence, every edge with label $>2$ (other than $[xy]$)
contains a linked vertex.

Now $\langle \delta'\rangle$ is conjugate to $\langle
\delta\rangle$ for some simplex $\delta$ of $V$. Then $\delta'$
contains more vertices than $\delta$. (If the odd edges of
$\delta'$ are collapsed to single vertices and each even $>2$
label of $V'$ is changed to 2, we obtain the (unique) diagram for
a right angled (all edge labels are 2) Coxeter group. The diagram
for this group is also obtained from $\delta$ if each even $>2$
label of $\delta$ is changed to 2. This latter description of this
diagram has the same number of vertices as $\delta$, but the
former diagram has fewer vertices (from the collapse of $[xy]$)
than $\delta '$.)

We obtain the desired contradiction by showing $\delta$ has at
least as many vertices as $\delta '$. Let $A'$ $(A)$ be the set
of vertices of $\delta '$ $(\delta)$ that belong to an edge of
$V'$ $(V)$ with label $>2$. As $\langle \delta'\rangle$ is
finite, no two adjacent edges of $\delta '$ have labels $>2$.
Similarly for $\delta$. The matching of Proposition \ref{P7} for
$\delta'$ and $\delta$ respects the matching for $V'$ and $V$.
Hence $\delta'$ contains an edge with label $>2$ iff $\delta $
contains the matching edge. So the number of vertices of $A'$
that belong to an edge of $\delta'$ with label $>2$ agrees with
the number of vertices of $A$ that belong to an edge of $\delta$
with label $>2$.

Suppose $[st]$ is an edge of $V'$ with label $>2$ and $s\in
\delta'$, $t\not \in \delta'$ and $s$ is not a vertex of an edge
in $\delta$ with label $>2$. Suppose $[ab]$ is the edge of $V$
matching $[st]$. Then $\{a,b\}\not \subset \delta$. Considering
the quotient of $W$ by $N(\delta)=N(\delta')$, we see
$\{a,b\}\cap \delta \not =\emptyset$. Hence we assume $b\in
\delta$ and $a\not \in \delta$. We wish to see that $b$ does not
belong to an edge of $\delta$ with label $>2$. Suppose $[bc]$ is
such an edge and $[uv]$ is an edge of $\delta'$ matching $[bc]$.
By Proposition \ref{P19} there is an odd edge path from $[st]$ to
$[uv]$. The first edge of this path cannot be $[tp]$, since then
$p\in \delta'$ and Lemma \ref{Tri} is violated. Hence the first
edge must be $[sp]$ and by assumption, $p\not\in \delta'$, so
$p\not \in \{u,v\}$. If $[pq]$ is the next edge, then $q\in
\delta'$ and again Lemma \ref{Tri} is violated. We conclude that
$b$ does not belong to an edge of $\delta$ with label $>2$.

Recall from Proposition \ref{P13} that if edges with label $>2$ of
$V'$ are adjacent, then their matching edges in $V$ are adjacent.
We show that $|A'|=|A|$ by verifying the following three
statements.

First, if $[st]$ and $[uv]$ are (non-adjacent) edges with labels
$>2$ of $V'$ such that $\{s,u\}\subset \delta'$, $t,v\not \in
\delta'$, and neither $s$ nor $u$ belongs to an edge of $\delta'$
with label $>2$, then the corresponding edges of $V$, call them
$[ab]$ and $[cd]$ respectively, are not adjacent. Otherwise,
there is an odd edge path from $[st]$ to $[uv]$. A contradiction
is obtained as in the former argument.

Suppose $[st]$ and $[sp]$ are edges of $V'$ with labels $>2$,
$s\in \delta'$, $\{t,p\}\subset V'-\delta'$ and $s$ not a vertex
of an edge of $\delta'$ with label $>2$. Then if $[ab]$ and
$[ac]$ are the edges of $V$ corresponding to $[st]$ and $[sp]$
respectively, either $a\in \delta$ and $a$ is not adjacent to an
edge of $\delta$ with label $>2$ or $\{b,c\}\subset \delta$ and
neither $b$ nor $c$ is adjacent to an edge of $\delta$ with label
$>2$. We show the latter scenario cannot occur. Otherwise, the
triangle $[abc]$ is contained in a maximal simplex $\tau$ of $V$,
and $\langle \tau\rangle$ is conjugate to $\langle \tau'\rangle$
for $\tau'$ a maximal simplex of $V'$. Hence $\{s,t,p\}$ forms a
triangle, no edge of which has odd label. By Lemmas 9 and 10,
there is no odd edge path from $\{t,p\}$ to $\delta'$ and no odd
edge path between $t$ and $p$. Let $K$ be the kernel of the map
of $\langle \sigma'\rangle$ to $\mathbb Z_2$ that takes
$\sigma'-\{s\}$ to 1 and $s$ to $-1$. A diagram for $\bar W\equiv
W/N(K)$ is obtained from $V'$ by removing the vertices
$\sigma'-\{s\}$ and all vertices that can be connected to
$\sigma'-\{s\}$ by an odd labeled edge path. The subgroup
$\langle s,t,p\rangle$ of $W$ injects under this quotient map. By
Lemma \ref{L29}, a diagram for $W/N(K)$ is obtained from $V$ by
removing some vertices of $\delta$ and identifying all others.
Neither $b$ nor $c$ are removed since Proposition \ref{P7}
(applied to $W/N(K)$ and the two diagrams for this group) implies
that the edges $[st]$ and $[sp]$ correspond to $[ab]$ and $[ac]$
respectively. Similarly $b$ and $c$ are not identified.

Finally, suppose $[st]$ and $[sp]$ are edges of $V'$ with labels
$>2$, $\{t,p\}\subset\delta'$, $s\in V'-\delta'$ and neither $t$
nor $p$ a vertex of an edge of $\delta'$ with label $>2$. Then if
$[ab]$ and $[ac]$ are the edges of $V$ corresponding to $[st]$
and $[sp]$ respectively, either $a\in \delta$ and $a$ is not
adjacent to an edge of $\delta$ with label $>2$ or
$\{b,c\}\subset \delta$ and neither $b$ nor $c$ is adjacent to an
edge of $\delta$ with label $>2$. We show the former scenario
cannot occur. As $\{s,t,p\}$ forms a triangle, no edge of this
triangle has odd label. If each odd labeled edge of $V'$ is
identified to a vertex then the resulting diagram is even with
Coxeter group $\bar W$ a quotient of $W$. Another diagram for
$\bar W$ is obtained from $V$ by changing labels on edges
corresponding to odd labeled edges to 2. The triangles $[stp]$ and
$[abc]$ induce triangles in the respective diagrams for $\bar W$
and the conjugate simplices groups $\langle \sigma'\rangle$ and
$\langle \sigma\rangle$ induce conjugate simplex groups in $\bar
W$. Since both of these diagrams are even, the previous argument
shows this is impossible.

By a completely analogous argument, we have:

\begin{lemma} \label{L31}
If $\sigma '\subset \delta'$ is such that $\langle \sigma
'\rangle$ is conjugate to $\langle \sigma\rangle$ for $\sigma
\subset \delta$, then $|A'\cap \sigma' |=|A\cap \sigma |$.
$\bullet$
\end{lemma}

Let $B'$ $(B)$ be the vertices of $\delta'$ $(\delta)$ that do
not belong to an edge with label $>2$. So $\delta '-A'=B'$ and
$\delta -A=B$. It suffices to show $|B'|\leq |B|$. If $b'\in B'$,
then there exists a simplex $\sigma_{b'}'\subset \delta'$ such
that $b'\in \sigma _{b'}'$, $\langle \sigma_{b'}'\rangle$ is
conjugate to $\langle \sigma _{b'}\rangle$ for some simplex
$\sigma _{b'}\subset \delta$ and not both $x$ and $y$ are in
$\sigma_{b'}'$. Since $\sigma _{b'}'$ can contain no linked
vertex, it must be right angled and so $|\sigma _{b'}'|=|\sigma
_{b'}|$. By Lemma \ref{L31}, $|\sigma _{b'}'\cap A'|=|\sigma
_{b'}\cap A|$, so $|\sigma_{b'}'\cap B'|=|\sigma_{b'}\cap B|$.

\begin{lemma}\label{L32}
If $\sigma_1',\ldots ,\sigma _n'$ are subsets of $\delta'$, and
$\langle \sigma _i'\rangle$ is right angled and conjugate to
$\langle \sigma _i\rangle$ for $\sigma _i\subset \delta$ then
$|B'\cap (\cap _{i=1}^n \sigma _i')|=|B\cap (\cap _{i=1}^n \sigma
_i)|$.
\end{lemma}

\noindent {\bf Proof:} Since $V$ is an even diagram, $\langle \cap
_{i=1}^n \sigma _i'\rangle$ is conjugate to $\langle
\sigma\rangle$ for $\sigma \subset \cap _{i=1}^n \sigma _i$.
Hence $|B'\cap (\cap _{i=1}^n \sigma _i')|\leq |B\cap (\cap
_{i=1}^n \sigma _i)|$, and it remains to show the reverse
inequality. We present the case $n=2$. The general case is
completely analogous. Assume $\langle\sigma _1\cap \sigma
_2\rangle$ is conjugate to $\langle \bar \sigma _1'\rangle$ for
$\bar \sigma _1'\subset \sigma _1'$, and also conjugate to
$\langle \bar \sigma _2'\rangle$ for $\bar \sigma _2'\subset
\sigma_2'$. As $\langle \bar \sigma _1'\rangle$ is conjugate to
$\langle \bar \sigma _2'\rangle$, if $v\in \bar \sigma _1'$ then
there is an odd edge path from $v$ to some vertex of $\bar \sigma
_2'$. But if $v\in B'$, it belongs only to edges labeled 2. Hence
$\bar \sigma _1'\cap B'=\bar \sigma _2'\cap B'=\bar \sigma
_1'\cap \bar \sigma _2'\cap B'$. Also, $|\sigma _1\cap \sigma
_2\cap B|=|\bar \sigma _i'\cap B'|=|\bar \sigma _1'\cap \bar
\sigma _2'\cap B'|\leq |\sigma _1'\cap \sigma _2'\cap B'|$.
$\bullet$

\medskip

The sets $\sigma _{b'}'\cap B'$ for $b'\in \delta '\cap B'$ cover
$B'$. (although it is not clear if the sets $\sigma_{b'}$ cover
$B$). Lemma \ref{L32} and the Inclusion-Exclusion Principle imply
$|B'|\leq |B|$ and the proof of Proposition \ref{P5} is complete.
$\bullet$

\section{Loops of Size 4}

\begin{proposition}\label{P33}
Suppose $(W,S)$ is an even Coxeter system with diagram $V$ and
$V'$ is another diagram for $W$. If $[xy]$, $[yz]$, $[zv]$ and
$[vx]$ are distinct edges of $V'$, and  $[xy]$ and $[yz]$ have odd
labels  then there is an edge (labeled 2) between $v$ and $y$.
\end{proposition}

\noindent {\bf Proof:} The edges $[vx]$ and $[vz]$ have label 2,
by Lemmas \ref{Tri} and \ref{Quo}. By Lemma \ref{Tri} there is no
edge $[xz]$. Suppose $V'$ is a minimal counterexample (with
respect to $T(V')$). We prove a collection of lemmas.

\begin{lemma} \label{L34}
\item {(1)} Every even label of an edge of $V'$ is 2.

\item {(2)} If $u\not =y$ is adjacent to $x$ (resp. $z$)then $[ux]$
(resp. $[uz]$) has label 2.

\item {(3)} Every odd labeled edge of $V'$ is adjacent to $y$ or $v$.

\item {(4)} If $[yu]$ (resp. $[vu]$) has odd label then
$(uv)^2=1$ (resp. $(uy)^2=1$).
\end{lemma}

\noindent {\bf Proof:} Otherwise a quotient map leads to a smaller
counterexample. $\bullet$

\medskip

Let $\bar V$ be the full subcomplex of $V'$ with vertex set
$\{v\}$ union all vertices of the odd labeled edges. (In
particular, $\{x,y,z,v\}\subset \bar V$.)

\begin{lemma}\label{L35}
Suppose $\langle \sigma '\rangle$ is conjugate to $\langle \sigma
\rangle$ for $\sigma'$ and $\sigma $ non-trivial simplices of
$V'$ and $V$ respectively.

\item {(1)} The simplex $\sigma '$ contains a vertex of $\bar V$.

\item {(2)} If $\sigma '$ contains a vertex of $V'-\bar V$ then
$\sigma '$ contains two vertices of $\bar V$.
\end{lemma}

\noindent {\bf Proof:} If $\sigma'$ contains no vertex of $\bar
V$ then $W/N(\sigma ')$ is a smaller counterexample. If $\sigma
'$ contains a vertex of $V'-\bar V$ and $t$ is the only vertex of
$\sigma '\cap \bar V$, then by Lemma \ref{L29} $W/N(\sigma
'-\{t\})$ is a smaller counterexample. $\bullet$

\begin{lemma}\label{L36}
Suppose $[x'y]$ is an odd labeled edge of $V'$ (so $(x'v)$ has
order 2). Then there is a vertex $t\in \bar V$ such that $(tv)$
has odd order, and for every simplex $\sigma '$ of $V'$
containing $\{x',v\}$ and such that $\langle \sigma '\rangle$ is
conjugate to $\langle \sigma \rangle$ for $\sigma$ a simplex of
$V$, $t\in \sigma '$. (By Lemma \ref{Tri}, $\sigma '$ cannot
contain a vertex $t'\ne t$ such that $(t'v)$ has odd order. In
this sense, $t$ is unique.)
\end{lemma}

\noindent {\bf Proof:} Suppose $\sigma '$ is a simplex of $V'$
containing $\{x',v\}$ and such that $\langle \sigma '\rangle$  is
conjugate to $\langle \sigma\rangle$ for some simplex $\sigma$ of
$V$, then $x,y,z \not\in \sigma '$. There is no vertex $t\in
\sigma '$ such that $[yt]$ is an odd labeled edge by Lemma
\ref{Tri}. If there is no vertex $t \in \sigma '$ such that
$(tv)$ has odd order, then (by Lemma \ref{L29}) $W/N(\{x'v\}\cup
(\sigma '-\{x',v\}))$ is an even Coxeter group and the diagram for
this Coxeter group obtained from $V'$ by identifying $x'$ and $v$,
and removing the vertices of $\sigma '-\{x',v\}$, contains a
triangle that violates Lemma \ref{Tri}.

Now suppose that $\sigma _1'$ and $\sigma _2'$ are simplices of
$V'$ as above, and $t_i$ is a vertex of $\sigma _i'$ such that
$t_iv$ has odd order and $t_1\not =t_2$. Then
$\{x',v\}\subset\sigma _1'\cap \sigma _2'\equiv \sigma '$. But
then there is a $t_3\in \sigma '$ such that $(t_3v)$ has odd
order. As $t_3$ is in both $\sigma _1'$ and $\sigma _2'$, we have
a contradiction to Lemma \ref{Tri}. $\bullet$

\medskip

Note that in the previous lemma $ty$ has order 2.

By a completely analogous argument we have:

\begin{lemma}\label{L37}
Suppose $[vt']$ is an odd labeled edge of $V'$ (so $(t'y)$ has
order 2). Then there is a vertex $x'\in \bar V$ such that $(x'y)$
has odd order, and for every simplex every simplex $\sigma '$ of
$V'$ containing $\{t',y\}$ and such that $\langle \sigma '\rangle$
is conjugate to $\langle \sigma \rangle$ for $\sigma$ a simplex
of $V$, $x'\in \sigma '$. (By Lemma \ref{Tri}, $\sigma '$ cannot
contain a vertex $x''\ne x'$ such that $(x''y)$ has odd order. In
this sense, $x'$ is unique.) $\bullet$
\end{lemma}

Let $x\equiv x_1$, $x_2, \ldots , x_n=z$ be the vertices of $\bar
V$ such that $(x_iy)$ has odd order. Let $t_i$ be the vertex of
Lemma \ref{L36} for $\{x_i,v\}$, so that $(t_iv)$ has odd order
and $(t_ix_i)$ has order 2.

Next we show that $t_i\not =t_j$ for $i\not =j$. Otherwise,
consider simplicies $\sigma _1'$ and $\sigma _2'$ containing
triangles $[t_ix_iy]$ and $[t_ix_jy]$ respectively, such that
$\langle \sigma_i'\rangle$ is conjugate to $\langle \sigma_i
\rangle$ for some simplex $\sigma _i$ of $V$. But then Lemma
\ref{L37} implies $x_i =x_j$, which is nonsense.

Similarly there is no edge between $t_i$ and $x_j$ for $i\ne j$
(see Figure 4). Hence we have:

\begin{lemma} \label{L38}
The only edges of $\bar V$ are $[vx_i]$, $[vt_i]$, $[yx_i]$,
$[yt_i]$ and $[x_it_i]$. $\bullet$
\end{lemma}

\medskip

\centerline{\bf Figure 4}

\medskip

\begin{lemma} \label{L39}
If $s$ is a vertex of $V'-\bar V$, then there is no pair of edges
$e_1, e_2$ such that $e_1$ connects $s$ to a point of
$\{x_i,t_i\}$ and $e_2$ connects $s$ to a point of $\{x_j,t_j\}$
for $i\not =j$.
\end{lemma}

\noindent {\bf Proof:}  Otherwise let $\sigma _1'$ and $\sigma
_2'$ be maximal simplices containing $e_1$ and $e_2$
respectively. Now, $s\in \sigma _1'\cap \sigma_2'\equiv \sigma
'$. By Lemma \ref{L35}, $\sigma '$ contains two vertices of $\bar
V$ and so $t_l$ or $x_l$ is an element of $\sigma '$ for some
$l$. But this is impossible by Lemma \ref{L38} $\bullet$

\begin{lemma} \label{L40}
The set $\{v,y\}$ separates $x$ from $z$.
\end{lemma}

\noindent {\bf Proof:} Suppose there is an edge path
$[xs_1],[s_1s_2],\ldots ,[s_mz]$ that does not intersect
$\{v,y\}$. By Lemma \ref{L39}, there is a smallest integer
$1<i\leq m$ such that $s_i$ is not connected to $x$ or $t_1$ by
an edge. Note that $s_{i-1}\not \in \bar V$ and $s_{i-1}$ is
connected to $x$ or $t_1$ by an edge. Also, by Lemma \ref{L39}
(with $s_{i-1}$ in place of $s$)  $s_i\not\in \bar V$. Let $\sigma
'$ be a maximal simplex of $V'$ containing $\{s_{i-1},s_i\}$ By
Lemma \ref{L35} there is $j\not =1$ and vertex $u\in
\{x_j,t_j\}\cap \sigma '$. But this is impossible by Lemma
\ref{L39} applied to $s_{i-1}$. $\bullet$

\medskip

By Proposition \ref{P13} and Lemma \ref{L23} the odd labeled edges
of $V'$ containing $v$ correspond to edges of $V$ with common
vertex $a$ and the odd labeled edges of $V'$ containing $y$
correspond to edges of $V$ with common vertex $b$. The subcomplex
of $V$ composed of the edges with label $>2$ and adjacent to $a$,
and the subcomplex of $V'$ composed of the edges with label $>2$
and adjacent to $b$ have trivial intersection by Proposition
\ref{P19}.

Let $\sigma _1'$ be a maximal simplex of $V'$ containing the
triangle $[x,t_1,v]$ and let $\sigma _2'$ be a maximal simplex of
$V'$ containing the triangle $[x,t_1,y]$. Let $\sigma _1$ and
$\sigma _2$ be simplices of $V$ such that $\langle \sigma
_i'\rangle$ is conjugate to $\langle \sigma _i\rangle$. Now,
$\sigma _1$ contains an edge $[ac]$ (with label $>2$)
corresponding to $[t_1v]$ and no vertex $c'\not =c$ such that
$[ac']$ has label $>2$. Similarly, $\sigma _2$ contains an edge
$[bd]$ corresponding to $[xy]$ and no vertex $d'\not =d$ such that
$[bd']$ has label $>2$. By conjugation we assume that $\langle
\sigma _1\rangle =\langle \sigma _1'\rangle$ and say $\langle
\sigma _2'\rangle =w\langle \sigma _2\rangle w^{-1}$ for some
$w\in W$. We have $\{x,t_1\}\subset \langle \sigma '
\rangle\equiv \langle \sigma _1'\rangle \cap \langle \sigma
_2'\rangle =\langle \sigma _1\rangle\cap w\langle \sigma
_2\rangle w^{-1}=\langle \sigma \rangle $ for $\sigma \subset
\sigma _1$. Since $V$ is even, $\sigma \subset \sigma _2$.  If
$s\in \sigma$ is a vertex of an edge of $V$ with label $>2$, then
$s\in\{a,b,c,d\}$. Since $\sigma'$ (and hence $\sigma$) is right
angled, $a$ or $c$ is not in $\sigma$, and $b$ or $d$ is not in
$\sigma$.  Since $W/N(\sigma ')=W/N(\sigma)$ is right angled $a$
and $b$ are elements of $\sigma$. We now have that $\langle
a,b\rangle$ is isomorphic to $ \mathbb Z_2\times \mathbb Z_2$.

In a completely analogous manner we find a simplex $\tau'$ of
$V'$ containing $\{z,t_n\}$ such that $\langle \tau '\rangle$ is
conjugate to $\langle \tau\rangle$ where $\tau$ is a simplex of
$V$ containing $\{a,b\}$. Let $C$ be the component of
$V'-\{v,y\}$ containing $\{x,t_1\}$ and $D=\Lambda'-C$ (so
$\{z,t_n\}\subset D$). $W$ decomposes as the amalgamated product
$\langle C\cup \{v,y\}\rangle \ast _{\langle v,y\rangle}\langle
D\rangle$. Since $\langle x,t_1\rangle \sim \mathbb Z_2\times
\mathbb Z_2$ is not a subgroup of a conjugate of the edge group
$\langle v,y\rangle \sim \mathbb Z_2\ast \mathbb Z_2$, $\langle
x,t_1\rangle $ can not stabilize a vertex of \cal {T} (the
Bass-Serre tree for $\langle C\cup \{v,y\}\rangle \ast _{\langle
v,y\rangle}\langle D\rangle$), other than $\langle C\cup
\{v,y\}\rangle$. Since $\{x,t_1\}\subset \sigma$ and since
$\langle \sigma\rangle$ stabilizes some vertex of \cal {T},
$\langle \sigma \rangle$ stabilizes $\langle C\cup
\{v,y\}\rangle$. Similarly, $\langle z,t_n\rangle$ only
stabilizes the vertex $\langle D\rangle$ of \cal {T}. As some
conjugate $w\langle \tau\rangle w^{-1}$ contains $\langle
z,t_n\rangle$, $w\langle \tau \rangle w^{-1}$  stabilizes
$\langle D\rangle$ (equivalently $w\langle \tau\rangle
w^{-1}<\langle D\rangle$). Hence, $\langle \tau\rangle<
w^{-1}\langle D\rangle w$ and so $\langle \tau \rangle$
stabilizes the vertex $w^{-1}\langle D\rangle$ of \cal {T}. But
then $\langle a,b\rangle\sim \mathbb Z_2\times \mathbb Z_2$
stabilizes distinct vertices of \cal {T}. This is impossible as
$\mathbb Z_2\times \mathbb Z_2$ is not a subgroup of $\mathbb
Z_2\ast \mathbb Z_2$. $\bullet$

\begin{proposition} \label{P41}
Suppose $W$ is a finitely generated Coxeter group with even
diagram $V$ and diagram $V'$ which is not even. If in $V'$, $[xy]$
has odd label, $[xu]$ and $[yz]$ have even labels,  and $[uz]$ has
odd label, then there is an (diagonal) edge $[xz]$ or $[yu]$.
\end{proposition}

\noindent {\bf Proof:} Otherwise assume $V'$ is a minimal
counterexample to the Proposition. Then we may assume that every
even edge of $V'$ is labeled 2 and by Proposition \ref{P33} and
Lemma \ref{Quo} there are no odd labeled edges, other than $[xy]$
and $[uz]$. Let $\sigma '$ be a maximal simplex of $V'$
containing $\{u,x\}$. Then $z,y\not\in \sigma'$. By Lemma
\ref{L29}, $W/N(\{ux\}\cup (\sigma '-\{u,x\}))$ is an even
Coxeter group with diagram obtained from $V'$ by identifying $u$
and $x$, and removing the vertices of $\sigma'-\{u,x\}$. This
contradicts Lemma \ref{Tri}. $\bullet$

\begin{lemma} \label {LA} Suppose $\mathbb Z_2^n=\langle
a,b,s_3,\ldots ,s_n\rangle =\langle a,b,e_3,\ldots ,e_n\rangle
\equiv G$. Then there is a retraction $h:G\to \langle a,b\rangle$
such that $h(a)=a$, $h(b)=b$, $h(e_i)\in \{a,b\}$, and $h(s_i)\in
\{a,1\}$ for all $i$.\
\end{lemma}

\noindent {\bf Proof:} Let $M_{n,n}$ be the coefficient matrix
for $s_i$. I.e. $s_i=m_{i,1}a+m_{i,2}b+m_{i,3}e_3+\cdots
+m_{i,n}e_n$, where each $m_{i,j}\in \{0,1\}$.

The matrix $M$ is invertible. As the block $[m_{i,j}]$ for $j>2$
and $i>2$ is invertible, a sequence of elementary row operations
(only involving rows 3 through $n$) can be used to transform $M$
to the matrix $\bar M$ where $\bar m_{i,i}=1$ for all $i$, $\bar
m_{i,j}=0$ for $i<j$ and for $i>j>2$.

We now define $h$. If $\bar m_{i,1}=1$ and $ \bar m_{i,2}=1$ then
$h(e_i)=b$.

If $\bar m_{i,1}=1$ and $ \bar m_{i,2}=0$ then $h(e_i)=a$.

If $\bar m_{i,1}=0$ and $ \bar m_{i,2}=1$ then $h(e_i)=b$.

If $\bar m_{i,1}=0$ and $ \bar m_{i,2}=0$ then $h(e_i)=a$.

To understand the effect this has on $h(s_i)$, some notation is
helpful. In $\bar M$, if $\bar m_{i,1}=1$, then replace it by
$1_a$. If $\bar m_{i,2}=1$, then replace it by $1_b$. If
$h(e_i)=a$ (respectively $b$) replace $\bar m_{i,i}$ by $1_a$
(respectively $1_b$). In each row of $\bar M$ (except the second)
there are an even number of $1_b$-entries. Reversing the above row
operations to obtain $M$ from $\bar M$ we see that if $\bar
m_{j,j}=1_a$ (respectively $1_b$) then $m_{i,j}\in \{0,1_a\}$,
(respectively $m_{i,j}\in \{0,1_b\}$) and $h(e_j)=a$
(respectively $h(e_j)=b$). Furthermore, each row of $M$ (except
the second) contains an even number of $1_b$-entries. I.e.
$h(s_i)\in \{1,a\}$. $\bullet$

\begin{proposition} \label {P43}
Suppose $W$ is a finitely generated Coxeter group with even
diagram $V$ and diagram $V'$ which is not even. If $[xy]$ has odd
label in $V'$ and $(xyst)$ is a simple loop in $V'$, then there is
an (diagonal) edge $[xs]$ or $[yt]$.
\end{proposition}

\noindent {\bf Proof:} Assume $V'$ is a minimal counterexample to
the proposition. By Propositions \ref{P33} and \ref{P41}, $[xy]$
is the only odd labeled edge of $(xyst)$. Each even edge of $V'$
is labeled 2. By Lemma \ref{Quo}, any odd labeled edge not
containing $x$ or $y$ must contain a vertex of $\{s,t\}$ and this
edge must be connected to the diagonally opposite vertex of
$\{x,y\}$ by an edge labeled 2. By Proposition \ref{P33}, there
can be at most one odd labeled edge containing a given vertex of
$\{x,y,s,t\}$.

\medskip

\noindent {\bf Case 1.} The only odd edge of $V'$ is $[xy]$.

\medskip

By Lemma \ref{L29} every simplex $\sigma '$ of $V'$ such that
$\langle \sigma '\rangle$ is conjugate to $\langle \sigma
\rangle$ for $\sigma $ a simplex of $V$ is either $\{u\}$ for
$u\in \{x,y,s,t\}$ or contains two vertices of $\{x,y,s,t\}$, but
not three. By Lemma \ref{L30} there is a vertex $u\in V'$ such
that every simplex $\sigma '$ of $V'$ containing $u$ and such that
$\langle \sigma'\rangle$ is conjugate to $\langle \sigma\rangle$
for $\sigma $ a simplex of $V$ contains $\{x,y\}$, and every
simplex $\tau'$ of $V'$ containing $\{x,y\}$ such that $\langle
\tau' \rangle$ is conjugate to $\langle \tau\rangle$ for $\tau$ a
simplex of $V$ contains $u$.

We show that $\{x,y\}$ separates $u$ from $t$. Choose a shortest
edge path avoiding $\{x,y\}$ from $u$ to $t$. Let $v$ be the
first vertex of this path that is not adjacent to both $x$ and
$y$ and let $v'$ be the previous vertex. Choose a maximal simplex
$\sigma _1'$ containing $\{x,y,v'\}$ and a maximal simplex
$\sigma _2'$ containing $\{v,v'\}$. Let $\sigma '=\sigma_1'\cap
\sigma _2'$. Assume $\sigma_2'$ does not contain $x$. Then
$\sigma'$ does not contain $x$, $s$ or $t$, which is impossible.

Now as $\{x,y\}$ separates $V'$, \cite{MT} implies that $W$
visually decomposes (with respect to $V'$) as $\langle
A\rangle\ast _{\langle x,y\rangle}\langle B\rangle$, where $A\cup
B$ is the vertex set of $V'$ and $A$ and $B$ properly contain
$\{x,y\}$. By \cite{MT}, there is a visual (with respect to $V$)
decomposition $C\ast _ED$ of $W$ such that $E$ is a subgroup of a
conjugate of $\langle x,y\rangle$. This implies that $E$ is
either trivial or $\langle v\rangle $ for some vertex $v$ of $V$.
Then $W$ visually decomposes (with respect to $V'$) as $\langle
F\rangle\ast _{\langle H\rangle }\langle G\rangle$ where $H$ is
trivial or $\{ v'\}$ for some vertex $v'$ of $V'$, $F\cup G$ is
the vertex set of $V'$ and $H$ is a proper subset of $F$ and $G$.
As $\langle x,y,s,t\rangle$ is 1-ended, we must have
$\{x,y,s,t\}$ a subset of $F$ or $G$. Assume $\{x,y,s,t\}\subset
F$. Let $e$ be a vertex of $G-H$ and $\sigma'$ be a maximal
simplex of $V'$ containing $e$. Then $\sigma '\cap \{x,y,s,t\}$
is either $\emptyset$ or $v'$. This is impossible as $\sigma'$
must contain two vertices of $\{x,y,s,t\}$. $\bullet$

\medskip

Hence there must be an odd edge at $s$ or $t$.

\medskip

\noindent {\bf Case 2.} Assume there is an odd edge $[tv]$ but no
odd edge at $s$.

\medskip

There is an edge $[yv]$ labeled 2 by Lemma \ref{Quo} and an edge
$[vx]$ labeled 2 by Proposition \ref{P41}. There are no odd
labeled edges, other than $[xy]$ and $[vt]$. Every simplex
$\sigma '$ (with at least two vertices) such that $\langle \sigma
'\rangle$ is conjugate to $\langle \sigma\rangle $ for some
simplex $\sigma$ of $V$ must contain at least two vertices of
$\{x,y,s,t,v\}$. There may or may not be an edge $[sv]$ labeled
2. By Lemma \ref{L29} there is a vertex $u$ of $V'$ such that if
$\sigma '$ is a simplex of $V'$, $\langle \sigma '\rangle $ is
conjugate to $\langle \sigma\rangle$ for some simplex $\sigma $
of $V$ and $x,y\subset \sigma'$ then $u\in \sigma '$. Also if
$\alpha '$ is a simplex of $V'$, $\langle \alpha '\rangle $ is
conjugate to $\langle \alpha \rangle$ for some simplex $\alpha $
of $V$ and $u\in \alpha'$ then $x,y\in \alpha '$. Note that
$u\not =v$.

Let $[ab]$ and $[cd]$ be the edges of $V$ corresponding to (see
Proposition \ref{P7}) $[xy]$ and $[tv]$ respectively. Note that
$\{a,b\}\cap \{c,d\}=\emptyset$ by Proposition \ref{P19}.

\begin{lemma} \label{L44}
$\{x,y,v\} $ separates $u$ from $t$ (and $s$).
\end{lemma}

\noindent {\bf Proof:}  Otherwise, there are consecutive vertices
$u=u_0,\ldots ,u_n=t$ such that no $u_i\in \{x,y,v\}$. Assume
that $i$ is the first integer such that $u_i$ does not commute
with $x$ and $y$. Suppose $u_i$ does not commute with $x$ (the
case $u_i$ does not commute with $y$ is completely analogous). Let
$\sigma '=\sigma _1'\cap \sigma _2'$ where $\sigma _1'$ is a
maximal simplex containing $\{x,y, u_{i-1}\}$ and $\sigma _2'$ is
a maximal simplex containing $\{u_{i-1}, u_i\}$. Then $t,s,x\not
\in \sigma '$ and so $u_{i-1},y,v\in \sigma '$.

Observe that $\sigma '-\{y\}\subset lk_2(x,y)$.

Suppose $ \sigma $ is a simplex in $V$ such that $\langle \sigma
\rangle$ is conjugate to $\langle \sigma '\rangle$.  By
conjugation, we may assume that $\langle \sigma \rangle =\langle
\sigma '\rangle$. We may assume that $a,c\in \sigma $ and
$b,d\not \in \sigma$, since $\langle \sigma\rangle$ is abelian.
Let $N=N(\sigma -\{a,c\})$. Note that $y,v\not \in N$.

Now we show:

\medskip

\noindent $(\ast)$ If $p\in \langle \sigma'-\{y\}\rangle$, then
$yp\not \in N$.

\medskip

\noindent Suppose $yp\in N$. As $yp$ conjugates
$xy(=w(ab)^2w^{-1})$ to $yx$, $w^{-1}ypw$ conjugates $(ab)^2$ to
$(ba)^2$. Proposition \ref{P12} implies that $w^{-1}ypw=ef$ for
$e\in lk_2(a,b)$ and $f$ of odd length in $\langle a,b\rangle$. As
$yp\in \langle \sigma\rangle$, $yp=ag$ for $g\in \langle \sigma
-\{a\}\rangle$. (Else, $yp=g\in \langle \sigma -\{a\}\rangle$
implying $ef=w^{-1}gw$ implying (the odd length element of
$\langle a,b\rangle$) $f=ew^{-1}gw$. But $ew^{-1}gw\in
N(S-\{a,b\})$ so this is impossible.) But $ag\not \in N$ and
$(\ast)$ is proved.

\medskip

Let $q:W\to W/N$ be the quotient map. Note that $q\langle
\sigma'\rangle = q\langle \sigma \rangle = \langle
q(a),q(c)\rangle \equiv Z_2\times Z_2$. As $y$, $v$ and $yv$ are
not elements of $N$, $q\langle \sigma'\rangle=\langle
q(y),q(v)\rangle$. To finish the Lemma, we show that $W/N$ is a
smaller counterexample.

First we show that if $m\in \sigma '-\{y,v\}$ then $q(m)\in
\{1,q(v)\}$: If $q(m)\not =1$, then $q(m)\in \{q(y), q(v),
q(yv)\}$. If $q(m)=q(y)$ then $q(my)=1$ implying $my\in N$,
contrary to $(\ast)$. If $q(m)=q(yv)$, then $q(mvy)=1$ implying
$mvy\in N$, contrary to $(\ast)$. Hence $q(m)=q(v)$.

If $K$ is the kernel of the restriction of $q$ to $\langle
\sigma'\rangle$, then the normal closure of $K$ in $W$ is $N$. As
$K$ is generated by $\{m\in \sigma' : q(m)=1\}\cup \{mv : m\in
\sigma '$ and $q(m)=q(v)\}$ , a diagram for the even Coxeter
group $W/N$ is obtained from $V'$ by removing the vertices of
$(\sigma'-\{y,v\})\cap ker(q)$ and identifying the remaining
vertices of $\sigma '-\{y,v\}$ with $v$. $\bullet$

\medskip

We can now finish Case 2. By \cite{MT} there is a simplex $\sigma
$ in $V$ such that $\sigma$ separates $V$ and $\langle \sigma
\rangle$ is a subgroup of a conjugate of $\langle x,y,v\rangle$.
Note that the edge $[ab]$ corresponding to $[xy]$ is not in
$\sigma$. (Since otherwise, the pigeon-hole principle implies
$\sigma =\{a,b\}$ and $\langle a,b\rangle$ is conjugate to
$\langle x,y,v\rangle$. This is impossible by since
$\{v,t\}\subset N(x,y,z)$, but $\langle c,d\rangle$ injects under
the quotient $W\to W/N(a,b)$.) Hence $\langle \sigma\rangle $ is
abelian and is either isomorphic to $\mathbb Z_2$ or $\mathbb
Z_2\times \mathbb Z_2$. Again applying \cite{MT}, either an edge
(labeled 2) or vertex of $V'$ must separate $V'$ and the group
generated by the vertices of this separating set is a subgroup of
a conjugate of $\langle x,y,v\rangle$. The set $\{x,y,s,t,v\}$
generates a 1-ended group and so no vertex $w$ of $V'$ can
separate $V'$. (Otherwise, there is a component $C$ of $V'-\{w\}$
such that $\{x,y,s,t,v\}\subset C\cup \{w\}$. But if $k$ is a
vertex of a component $K\ne C$ of $V'-\{w\}$, then a maximal
simplex of $V'$ containing $k$ could contain at most one vertex of
$\{x,y,s,t,v\}$, which is impossible.) If an edge labeled 2
separates $V'$ and the group for this edge is conjugate to a
subgroup of $\langle x,y,v\rangle$, then the only candidates are
the edges $[xv]$, $[yv]$ and $[xt]$. The groups for these edges
are all conjugate. All cases are completely similar, and we
assume $[yv]$ separates $V'$. Then (by \cite{MT}) there is an edge
separating $V$ such that the corresponding group is conjugate to
$\langle y,v\rangle$. We may assume that this edge is $[ac]$. By
conjugation, we may assume that $\langle y,v\rangle=\langle
a,c\rangle$. Note that $vy$ conjugates $xy(=w_1(ab)^2w_1^{-1})$
to $yx(=w_1(ba)^2w_1^{-1})$.

If $vy=c$, then $cw_1(ab)^2w_1^{-1}c=w_1(ba)^2w_1^{-1}$. But in
$W/N(\{c\})$, $(ab)^2\ne(ba)^2$, so this is impossible. If
$vy=a$, then choose $g\in \langle x,y\rangle$ such that
$gvyg^{-1}=xv$. As $xv$ conjugates $tv(=w_2(cd)^2w_2^{-1})$ to
$vt(=w_2(dc)^2w_2^{-1})$, $gag^{-1}$ conjugates
$w_2(cd)^2w_2^{-1}$ to $w_2(dc)^2w_2^{-1}$. But in $W/N(\{a\}$,
$(cd)^2\ne (dc)^2$, so this is impossible. The only other
possibility is $vy=ac$. Hence we have $\{y,v\}=\{a,c\}$.

The group $\langle x,y,s,t,v\rangle$ is 1-ended and hence
$\{x,y,s,t,v\}$ is a subset of $\{v,y\}$ union a component $C$ of
$V'-\{v,y\}$. Let $\sigma'$ be a maximal simplex containing a
vertex of a component of $V'-\{v,y\}$ other than $C$. Then
$\sigma'$ must also contain $\{v,y\}(=\{a,c\})$ and each edge
label of $\sigma'$ is 2. Assume that $\langle \sigma '\rangle
=w\langle \sigma \rangle w^{-1}$ for $\sigma$ a simplex of $V$.
Then $a,c\in \sigma$, $waw^{-1}=a$ and $wcw^{-1}=c$.  We wish to
apply Lemma \ref{LA} to obtain a smaller counterexample.

Write $\langle \sigma '\rangle =\langle a,c,e_3,\ldots ,
e_n\rangle=\langle a,c,ws_3w^{-1},\ldots
,ws_nw^{-1}\rangle=\langle w\sigma w^{-1}\rangle$  and suppose $h$
is the retraction of this group to $\langle a,c\rangle$ defined in
Lemma \ref{LA}. If $h(ws_iw^{-1})=a$, then observe that
$ws_iw^{-1}a=w(s_ia)w^{-1}$. Thus, an even diagram for
$W/N(ker(h))$ is obtained from $V$ by removing all $s_i$ such
that $h(ws_iw^{-1})=1$ and identifying $s_i$ with $a$ if
$h(ws_iw^{-1})=a$.  Another diagram for $W/N(ker(h))$ is obtained
from $V'$ by identifying $e_i$ with $v$ when $h(e_i)=v$ and
identifying $e_j$ with $y$ when $h(e_j)=y$. Then $W/N(ker(h))$ is
a smaller counterexample, finishing Case 2.

Now the final case.

\medskip

\noindent {\bf Case 3.} Assume there are edges $[tv]$ and $[su]$
with odd labels. Then there are edges $[uy]$ $[ux]$ $[vy]$ and
$[vx]$ with labeled 2. (See Figure 5)

\medskip

\centerline{\bf Figure 5}

\medskip

Let $\{a,b,c,d,e,f\}\subset V$ be vertices such that the edge
correspondence of Proposition \ref{P7} relates $[xy]$ to $[ab]$,
$[tv]$ to $[cd]$ and $[us]$ to $[ef]$. By Proposition \ref{P19},
the edges $[ab]$, $[cd]$ and $[ef]$ are mutually disjoint.  By
Theorem \ref{T26} there is an edge $[gh]\in \{[ab],[cd],[ef]\}$
such that every other edge of $V$ containing $g$ is labeled 2 and
such that if $[gk]$ is such an edge, then $[kh]$ is also an edge
of $V$. We call $g$ a {\it special} vertex of $[gh]$. A quotient
argument shows that each of $[ab],[cd],[ef]$ contains a special
vertex. E.g. if the edge labels on $[cd]$ and $[ef]$ are changed
to 2, then Theorem \ref{T26} implies that $[ab]$ must have a
special vertex. The next lemma implies that there cannot be an
edge $[uv]$ in the minimal counterexample.

\begin{lemma} \label{L45}
Suppose $(W,S)$ is a finitely generated even Coxeter
system,  $V'$ is a  diagram for $W$ with non-intersection odd
edges $[xy]$, $[tv]$ and $[us]$, and even edges $[xu]$, $[xv]$,
$[xt]$ $[yu]$, $[ys]$, $[yv]$, $[ts]$ and $[uv]$, (and so a
tetrahedron $[xyuv]$ and triangles $[xtv]$ and $[suy]$). Then
there is an edge (labeled 2) $[yt]$ or $[xs]$.
\end{lemma}

\noindent {\bf Proof:} If an edge not listed in the hypothesis,
between two vertices of $\{x,y,u,v,t,s\}$ exists, it must have
label 2 by Lemma \ref{Tri}. Assume $V'$ is a minimal
counterexample to the Lemma. All even edges of $V'$ are labeled 2
by the minimality of $V'$. Every odd edge of $V'$ contains a
vertex of $\{x,y,u,v,t,s\}$ (otherwise collapse for a smaller
counterexample). Suppose $V'$ contained an odd edge other than
$[xy]$, $[tv]$ or $[us]$. Then this edge must contain $x$, $y$,
$s$ or $t$ (otherwise collapse).  If $[xw]$ is an odd edge for
$w\not =y$, there must be an edge $[ws]$ (or collapsing would
give a smaller counterexample). In this situation, Proposition
\ref{P33} implies there is an edge $[xs]$ and we are finished.
Similarly if there is an odd edge at $y$, $s$ or $t$. Hence we
may assume there is no odd edge of $V'$ other than $[xy]$, $[tv]$
and $[us]$.

We assume $a$, $d$ and $f$ are special vertices of $[ab]$, $[cd]$
and $[ef]$ respectively. Observe that if $\sigma$ is a simplex of
$V$ containing $a$ then $\sigma \cup \{b\}$ is a simplex of $V$.
Similarly for $d$ and $f$. Let $\sigma _1'$ be a maximal simplex
of $V'$ containing $\{x,y,v,u\}$. Then $\langle \sigma _1'\rangle$
is conjugate to $\langle \sigma_1 \rangle$ for $\sigma_1 $ a
(maximal) simplex of $V$. Now $a,b\in \sigma_1$ and $c,e\in
\sigma_1$ (as $\{a,b\}$ does not commute with $\{c,d\}$ or
$\{e,f\}$). Similarly, considering the triangles $[yus]$ and
$[tvx]$, we see that $\{b\}$ commutes with $\{e,f,c,d\}$.

Consider the loop $(tvus)$. By Proposition \ref{P41}, there must
be an edge $[tu]$ or $[vs]$. If both exist, then $\{c,d\}$
commutes with $\{e,f\}$ and there is a maximal simplex of $V$
containing $\{b,c,d,e,f\}$, implying there is a maximal simplex
in $V'$ containing $\{t,v,u,s\}$ and either $x$ or $y$, which is
impossible in our minimal counterexample.

Now suppose $[vs]$ is an edge of $V'$, but $[tu]$ is not. (See
Figure 6.)

Let $\sigma _2'$ be a maximal simplex of $V'$ containing
$\{u,s,y,v\}$. Then $\langle \sigma _2'\rangle$ is conjugate to
$\langle \sigma_2 \rangle$ for $\sigma_2 $ a (maximal) simplex of
$V$ and $\{e,f,b,c\}\subset \sigma _2$. Choose a maximal simplex
$\sigma _3'$ of $V'$ containing $\{t,v,x\}$. Then there exists
$\sigma _3$, a maximal simplex of $V$ containing $\{c,d,b\}$ and
such that  $\langle \sigma _3'\rangle$ is conjugate to $\langle
\sigma_3 \rangle$. Let $\sigma _4'$ be a maximal simplex of $V'$
containing $\{t,v,s\}$. Then there exists $\sigma _4$, a maximal
simplex of $V$ containing $\{c,d\}$ and either $f$ or $e$, and
such that $\langle \sigma _4'\rangle$ is conjugate to $\langle
\sigma_4 \rangle$. But then either $\{b,c,d,f\}$ or $\{b,c,d,e\}$
is a simplex. This implies there is a simplex $\sigma'$ of $V'$
containing $\{t,v\}$, a vertex of $\{x,y\}$ and a vertex of
$\{u,s\}$. This is impossible and $[vs]$ is not an edge of $V'$.

\medskip

\centerline{\bf Figure 6}

\medskip

Next assume that $[tu]$ is an edge, but $[sv]$ is not an edge of
$V'$. (See Figure 7.)

We may assume that $\{x,y,u,v\}\subset \sigma _1'$ and
$\{a,b,c,e\}\subset \sigma _1$ as above. Also,
$\{x,u,v,t\}\subset \sigma _2'$ and $\{c,d,b,e\}\subset \sigma
_2$; $\{u,s,y\}\subset \sigma _3'$ and $\{e,f,b\}\subset \sigma
_3$; and $\{t,u,s\}\subset \sigma _4'$ and $\sigma _4$ contains
$\{e,f\}$ and either $c$ or $d$. But then either $\{e,f,c,b\}$ or
$\{e,f,d,b\}$ is a simplex.

\medskip

\centerline{\bf Figure 7}

\medskip

This implies there exists a simplex $\sigma '$ of $V'$ containing
$\{u,s\}$, a vertex of $\{x,y\}$ and a vertex of $\{v,t\}$. This
is impossible. $\bullet$

\medskip

\noindent {\bf Remark 4} Consider the loop $(suxt)$. There must be
an (labeled 2) edge $[ut]$ or $[uv]$, or collapsing $[vt]$ gives a
smaller counterexample. Consider the loop $(syvt)$. There must be
an (labeled 2) edge $[vs]$ or $[uv]$, or collapsing $[us]$ gives a
smaller counterexample. Either of these observations (along with
Lemma \ref{L45})  could be used in conjunction with the ideas of
\S 7 to complete a proof of Case 3.

\medskip
Lemma \ref{L45} implies $[uv]$ is not an edge of $V'$. At this
point Figure 8 is our model.

\medskip

\centerline{\bf Figure 8}

\medskip

Again if $\sigma '$ is a simplex of $V'$ such that $\langle
\sigma '\rangle$ is conjugate to a $\langle \sigma\rangle$, for
$\sigma$ a simplex of $V$, and $\sigma'$ contains a vertex of
$V'-\{x,y,t,s,u,v\}$, then $\sigma'$ must contain two vertices of
$\{x,y,t,s,u,v\}$.

\begin{lemma} \label{L46}
The set $\{x,y,u,v\}$ separates $V'$.
\end{lemma}

\noindent {\bf Proof:} Suppose not. By Lemma \ref{L30}, there is
a vertex $z$ of $V'-\{x,y\}$ such that every simplex $\sigma '$
of $V'$ containing $[xy]$ and such that $\langle \sigma '\rangle$
is conjugate to $\langle \sigma \rangle$ for $\sigma$ a simplex of
$V$, contains $z$ and every simplex $\sigma '$ of $V'$ containing
$z$ such that $\langle \sigma '\rangle$ is conjugate to $\langle
\sigma \rangle$ for $\sigma$ a simplex of $V$, contains $[xy]$.
Clearly $z\not\in \{v,u,t,s\}$.  Suppose $z=z_0,z_1,\ldots ,
z_n=t$ are consecutive vertices of $V'-\{x,y,u,v\}$. Assume that
$i$ is the first integer such that $z_i$ is not adjacent to both
$x$ and $y$. Both cases are analogous, and we assume that $z_i$
is not adjacent to $y$. Let $\sigma_1'$ be a maximal simplex of
$V'$ containing $\{x,y,z_{i-1}\}$, $\sigma_2'$ a maximal simplex
containing $\{z_{i-1}, z_i\}$. Let $\sigma'=\sigma _1'\cap \sigma
_2'$. Then $\{y,s,t\} \cap\sigma '=\emptyset$. Then $\sigma '$
must contain $z_{i-1}$ and either $\{x,u\}$ or $\{x,v\}$  ($[vu]$
is not a possibility by Lemma \ref{L45}). Both cases have
analogous proofs and we assume $\{x,u\}\subset \sigma '$. Say
$\sigma $ is a simplex of $V$ and $\langle \sigma \rangle$ is
conjugate to $\langle \sigma '\rangle$. Then we may assume
$b,e\in \sigma $ and $a,c,d,f\not \in \sigma$. Let $N=N(\sigma
-\{b,e\})$. Each vertex of $\sigma '-\{x\}$ commutes with $x$ and
$y$. An argument completely analogous to that for statement
$(\ast)$ in the proof of Lemma \ref{L44} implies that if $p\in
\langle \sigma '-\{x\}\rangle$, then $xp\not\in N$. Just as in the
argument following the proof of $(\ast)$, this implies $W/N$ is a
smaller counterexample. $\bullet$

\medskip

By \cite{MT} there is a full subgraph $A$ separating $V$ such that
$\langle A\rangle$ is conjugate to a subgroup of $\langle
u,v,x,y\rangle =\langle u,x,y\rangle\ast _{\langle x,y\rangle}
\langle u,x,y\rangle$.  The edge $[ab]$ is not in $A$, for
otherwise, a conjugate of $\langle a,b\rangle$ is a subgroup of
$\langle u,x,y\rangle$ or $\langle v,x,y\rangle$. But all three
of these groups have the same order. This would imply $\langle
a,b\rangle$ is conjugate to $\langle u,x,y\rangle$ or $\langle
v,x,y\rangle$, but clearly $N(\{a,b\})$ is not equal to
$N(\{u,x,y\})$ or $N(\{v,x,y\})$. So, $A$ is right angled. The
group $\langle x,y,u,v\rangle = \langle x,y,v\rangle \ast
_{\langle x,y\rangle }\langle x,y,u\rangle$ is 2-ended and
contains no copy of $\mathbb Z_2\times \mathbb Z_2\times \mathbb
Z_2$. This implies that $\langle A\rangle= 1$, $\mathbb Z_2$,
$\mathbb Z_2\times \mathbb Z_2$, $\mathbb Z_2 \ast \mathbb Z_2$,
or $\mathbb Z_2\times \mathbb Z_2\ast_{\mathbb Z_2}\mathbb
Z_2\times \mathbb Z_2$. Using \cite{MT} again, a subset $B$ of the
vertices of $V'$ separates $V'$ such that $\langle B\rangle$ is
right angled and a subgroup of a conjugate of $\langle A\rangle $
and so a subgroup of a conjugate of $\langle x,y,u,v\rangle$.
This implies that $B$ is a subset of $\{u,x,v\}$, $\{u,y,v\}$,
$\{ x,t \}$, or $\{ y,s \}$. By Corollary 6 of \cite{MT} and we
may assume that $\langle A\rangle$ is conjugate to $\langle
B\rangle$. Note that $\{x,y,u,v,s,t\}$ is a subset of $B$ union a
component of the compliment of $B$. Select a maximal simplex
$\sigma '$ intersecting another component of the compliment of
$B$. Then $\sigma '\cap B$ must be equal to $\{x,t\}$, $\{x,v\}$,
$\{x,u\}$, $\{y,v\}$, $\{y,s\}$ or $\{y,u\}$ and be conjugate to
a simplex of $V$. Now proceed as in Case 2. $\bullet$

\section{The Proof of Theorem 2}

In this section we finish the proof of Theorem \ref{T2}.
Throughout this section, we assume that $V'$ is a smallest
diagram for an even Coxeter group such that $V'$ contains a
simple loop $l$ without shortcuts, the length of $l$ is $\geq 5$
and $l$ contains an odd labeled edge $[xy]$. By the minimality of
$V'$, all even edges of $V'$ are labeled 2. By Lemma \ref{Quo},
Proposition \ref{P43} and the minimality of $V'$, all edges of
$l$ other than $[xy]$ are labeled 2.

\begin{lemma} \label{L47} If $l'$ is a loop of $V'$ containing 2
odd labeled edges, then $l'$ must have a shortcut.
\end{lemma}

\noindent {\bf Proof:} Otherwise, the diagram obtained from $V'$
by collapsing one of the edges of $l'$ contradicts Proposition
\ref{P43} or is a smaller example than $V'$. $\bullet$

\begin{lemma} \label{L48}
If $[xu]$ has odd label, then $u=y$.
\end{lemma}

\noindent {\bf Proof:} By Lemma \ref{Quo}, $u$ is connected to a
vertex of $l-st(x)$ by an edge. Say the consecutive vertices of
$l$ are $x=a_0, a_1, \ldots , a_n=y$. Let $i$ be the largest
integer such that $[ua_i]$ is an edge of $V'$. By Lemma \ref{Tri}
and Proposition \ref{P33}, $i<n-1$. The loop with consecutive
vertices $(yxua_i\ldots a_{n-1})$ has no shortcuts, contradicting
Lemma \ref{L47}. $\bullet$

\begin{lemma} \label{L49}
If $[uv]\ne [xy]$ is an odd labeled edge, then $[uv]$ has one
vertex on $l$ and the other vertex in $lk_2(x,y)$.
\end{lemma}

\noindent {\bf Proof:} If neither $u$ nor $v$ is a vertex of $l$,
then $W/N(uv)$ is a smaller example. Assume $v$ is a vertex of
$l$. If $v$ is not adjacent to $x$ or $y$ (and $u\not\in
lk_2(x,y)$), then the quotient of $W$ by $N(uv)$ gives a smaller
example. If say $v$ is adjacent to $x$, then $u$ must be adjacent
to $y$ or again $W/N(uv)$ is a smaller example. Now by Lemma
\ref{L47} (applied to the loop $(vuyx)$), $u$ is adjacent to $x$.
$\bullet$

\begin{lemma} \label{L50}
If $[uv]$ is an odd labeled edge of $V'$ then $[uv]$ is contained
in a simple loop of length $\geq 5$ without shortcuts.
\end{lemma}

\noindent {\bf Proof:} Suppose otherwise. We assume that $v\in l$
and $u\in lk_2(x,y)$.  If $s,t$ are the vertices of $l$ adjacent
to $v$, then $s,t\in lk_2(u)$ or $[uv]$ belongs to a loop of
length $\geq 5$ without shortcuts (all edges of this path not
containing $u$ would be in $l$).  Note that no vertex of
$l-\{s,t\}$ belongs to $lk_2(v)$ (and so no such vertex belongs
to $lk_2(u,v)$).

Next we show that if $[vw]$ is such that $w\ne u$ and $w\not\in
lk_2(u,v)$, then there is no edge path in $V'-lk_2(u,v)$ from $w$
to a vertex of $l$. Otherwise, there is a simple edge path loop
containing $[uv]$ and avoiding $lk_2(u,v)$. A shortest such loop
contradicts the hypothesis on $[uv]$.

Now, let $U$ be the union of all components $K$ of $V'-lk_2(u,v)$
such that for some vertex $w\in K$, $[wv]$ is an edge. We have
$U\cap l=\emptyset$. Twist $U$ around $[uv]$ to form the diagram
$V''$ for $W$. If $z\not =u$ and $[vz]$ is an edge of $V''$, then
$z\in lk_2(u,v)$. By Proposition \ref{P5} there is a vertex $w$
($\not\in \{s,t\})$ of $V''$ such that the triangle $[uvw]$ can be
replaced by an edge $[uz]$ with even label $>2$. The resulting
diagram for $W$ is smaller than $V'$ and retains $l$ (with $v$
replaced by $z$). $\bullet$

\begin{lemma} \label{L51}
If $v$ is a vertex of $l$ not adjacent to $x$ or $y$, then there
is no odd edge at $v$.
\end{lemma}

\noindent {\bf Proof:} Assume that $[uv]$ is such an odd edge. By
Lemma \ref{L49}, $u\in lk_2(x,y)$. Let $l'$ be a simple edge path
without shortcuts containing $[uv]$ and having length $\geq 5$.
Then $x$ or $y$ is a vertex of $l'$ (otherwise, $W/N(xy)$ is a
smaller example). We may assume $x$ is a vertex of $l'$. By Lemma
\ref{L49}, $y\in lk_2 (u,v)$. This is impossible as $[yv]$ would
be a shortcut in $l$. $\bullet$

\begin{lemma} \label{L52}
If $s\ne y$ and $t\ne x$ are vertices of $l$ adjacent to $x$ and
$y$ respectively, then there is not an odd labeled edge at $s$ and
an odd labeled edge at $t$.
\end{lemma}

\noindent {\bf Proof:}  Otherwise, say $[su]$ and $[tv]$ have odd
labels. By Lemmas \ref{L50} and \ref{L48} $u\ne v$. By Lemma
\ref{L49}, $u,v\in lk_2(x,y)$. Assume $l_s$ and $l_t$ are simple
loops without shortcuts and of length $\geq 5$ containing $[su]$
and $[tv]$ respectively.  By Lemma \ref{L49}, $l_s$ and $l_t$
contain $y$ and $x$ respectively. Hence the paths $(suy)$ and
$(tvx)$ are contained in $l_s$ and $l_t$ respectively. Since
there is no edge from $t$ to $s$, Lemma \ref{L49} implies that
$v\in lk_2 (s,u)$. By Lemma \ref{L49}, $l_t$ contains $s$ or $u$
and so the vertex of $l_t$ following the path $(tvx)$ must be $s$
or $u$, but this is impossible as both are connected to $v$ by an
edge (creating a shortcut in $l_t$.) $\bullet$

\medskip

We can now complete the reduction.

\medskip

\noindent {\bf Case 1.} Suppose the only odd labeled edge of $V'$
is $[xy]$.

\medskip

Then say $(xyst)$ is a subpath of $l$. Let $\sigma'$ be a simplex
of $V'$ containing $\{s,t\}$ such that $\langle \sigma'\rangle $
is conjugate to $\langle \sigma \rangle$ for some simplex $\sigma
$ of $V$. By Lemma \ref{L29}, the quotient of $W$ by $N(\{st\}\cup
(\sigma'-\{s,t\}))$ is an even Coxeter group. A diagram for this
group is obtained from $V'$ by identifying $s$ and $t$ and
removing the vertices of $\sigma '-\{s,t\}$. This diagram
contains a loop of length $\geq 4$, with odd labeled edge $[xy]$
and no shortcuts. This diagram is smaller than $V'$,
contradicting Proposition \ref{P43} or the minimality of $V'$.

\medskip

\noindent {\bf Case 2.} Suppose $V'$ contains exactly two odd
labeled edges.

\medskip

Say $(uxyst)$ is a subpath of $l$ and $[uv]$ is an odd labeled
edge. By Lemma \ref{L49}, $v\in lk_2(x,y)$. If $[cd]$ is an edge
of $l$ such that $\{c,d\}\cap \{u,x,y\}=\emptyset$, and $\sigma'$
is a simplex of $V'$ containing $[cd]$ such that $\langle \sigma
'\rangle$ is conjugate to $\langle \sigma\rangle$ for some
simplex $\sigma $ of $V$, then $\sigma'$ must contain $v$
(otherwise, the quotient of $W$ by $N(\{cd\}\cup
(\sigma'-\{c,d\})$ gives a smaller example). In particular, $v$
is connected to each vertex of $l$ by an edge.

Let $l'$ be a simple edge loop in $V'$ of length $\geq 5$,
containing $[uv]$, and without shortcuts. Then $l\cap
l'=\{u,y\}$. Note that $(uvy)$ is a subpath of $l'$. Let $\sigma'$
be a simplex of $V'$ containing $[stv]$ such that $\langle \sigma
'\rangle$ is conjugate to $\langle \sigma\rangle$ for some
simplex $\sigma $ of $V$. Then $u,x,y\not\in \sigma'$ and the
only vertex of $l'$ in $\sigma '$ is $v$. A diagram for the even
Coxeter group $W/N(\sigma'-\{v\})$ is obtained from $V'$ by
removing the vertices of $\sigma '-\{v\}$. But this diagram
contains a faithful copy of $l'$ and so contradicts the minimality
of $V'$. The proof of Theorem \ref{T2} is complete.

\end{document}